\documentclass[preprint,fleqn,sort&compress,3p]{elsarticle}
\usepackage{mathrsfs}
\usepackage{amsfonts}
\usepackage{graphicx}
\usepackage{url}
\usepackage{amssymb}
 \usepackage{amsthm,amsmath}
\usepackage{bm}
\usepackage{latexsym}
\usepackage{booktabs}
\usepackage{hyperref}
\usepackage{breakurl}
\usepackage{algorithm,algorithmic}

\numberwithin{equation}{section}
\newtheorem{theorem}{Theorem}[section]
\newtheorem{proposition}{Proposition}[section]

\newtheorem{lemma}{Lemma}[section]

\journal{arXiv preprint}
\usepackage{exscale}
\usepackage{relsize}
\usepackage{color}
\usepackage{rotating}

\begin{document}
\title{A fast implicit difference scheme for solving the generalized time-space fractional
diffusion equations with variable coefficients}


\author[a,b]{Xian-Ming Gu}
\ead{guxianming@live.cn, guxm@swufe.edu.cn}
\author[c]{Ting-Zhu Huang\corref{cor1}}
\ead{tingzhuhuang@126.com}
\cortext[cor1]{Corresponding author}
\author[c]{Yong-Liang Zhao}
\ead{ylzhaofde@sina.com}
\author[a]{Pin Lyu}
\ead{plyu@swufe.edu.cn}
\author[e]{Bruno Carpentieri}
\ead{bcarpentieri@gmail.com}
\address[a]{School of Economic Mathematics/Institute of Mathematics, \\ Southwestern University of Finance
and Economics, Chengdu 611130, Sichuan, P.R. China}
\address[b]{Bernoulli Institute of Mathematics, Computer Science and Artificial Intelligence, \\
University of Groningen, Nijenborgh 9, P.O. Box 407, 9700 AK Groningen, The Netherlands}
\address[c]{School of Mathematical Sciences, \\ University of Electronic
Science and Technology of China, Chengdu, Sichuan 611731, P.R. China}
\address[e]{Facolt\`{a} di Scienze e Tecnologie informatiche,\\
Libera Universit\`{a} di Bolzano, Dominikanerplatz 3 - piazza Domenicani, 3 Italy - 39100, Bozen-Bolzano}

\begin{abstract}
In this paper, we first propose an unconditionally stable implicit difference scheme for solving
generalized time-space fractional diffusion equations (GTSFDEs) with variable coefficients. The
numerical scheme utilizes the $L1$-type formula for the generalized Caputo fractional derivative
in time discretization and the second-order weighted and shifted Gr\"{u}nwald difference (WSGD)
formula in spatial discretization, respectively. Theoretical results and numerical tests are conducted to verify
the $(2 - \gamma)$-order and 2-order of temporal and spatial convergence with $\gamma\in(0,1)$ the
order of Caputo fractional derivative, respectively. The fast sum-of-exponential approximation of
the generalized Caputo fractional derivative and Toeplitz-like coefficient matrices are also developed
to accelerate the proposed implicit difference scheme. Numerical experiments show the effectiveness
of the proposed numerical scheme and its good potential for large-scale simulation of GTSFDEs.
\end{abstract}
\maketitle
%
\section{Introduction}
\label{sec1}
In recent years there has been a growing interest in the field of fractional calculus. For instance, Podlubny
\cite{Podlubny99}, Samko \emph{et al.} \cite{Samko93} and Kilbas \emph{et al.} \cite{Kilbas06} provide the
history and a comprehensive treatment of this subject. Many phenomena in engineering, physics, chemistry and
other sciences can be described very successfully by using fractional partial differential equations (FPDEs).
Diffusion with an additional velocity field and diffusion under the influence of a constant external force field
are, in the Brownian case, both modelled by the diffusion equation. In the case of anomalous diffusion this is
no longer true, i.e., the space fractional generalization may be different for the transport in external force
field \cite{Metzler}. Under the framework of the continuous time random walks (CTRWs) model, the fractional diffusion,
Fokker-Planck and Feynman-Kac equations \cite{Metzler,Xu2016} can be derived with power-law waiting time distribution
(WTD), assuming the particles may exhibit long waiting time. However, for some practical physical processes, it is
necessary to make the first moment of the waiting time measure finite. This leads to the generalized time fractional
diffusion equation corresponding to the CTRWs model with some more complicated WTDs (beyond the power-law limit)
\cite{Sandev15,Alikhanov17,Santos19}, e.g., the tempered \cite{Hanert14,Sabzikar15,Zhang17,Chen2018,Guo2019} and
the scale-weight \cite{Xu2013,Gao2017} power law WTDs. In one word, the generalization of time-space fractional
diffusion equations where the sub-diffusion in time and the super-diffusion in space simultaneously \cite{Chi2017}
will be meaningful to model the anomalous diffusion with complicated physical processes.

Based on the above considerations, in this work, we are interested in developing fast numerical methods for solving the initial-boundary value
problem of the generalized time-space fractional diffusion equation (GTSFDE) with variable coefficients
\begin{equation}
\begin{cases}
{}^{C}_{0}D^{\gamma,\lambda(t)}_{t}u(x,t) = \xi(x,t)\Big[p{}_{x_L}D^{\alpha}_{x}u(x,t)
+ (1 - p){}_xD^{\alpha}_{x_R}u(x,t)\Big] + f(x,t),& (x,t)\in(x_L,x_R)\times(0,T),\\
u(x,0) = \phi(x), & x\in [x_L,x_R],\\
u(x_L,t) = \varphi(t),\quad u(x_R,t) = \psi(t),& t\in(0,T],\\
\end{cases}
\label{eq1.1}
\end{equation}
where $\alpha\in(1,2]$, $\gamma\in(0,1)$. The parameter $p\in[0,1]$, called ``skewness", represents
the proportion of high-velocity ``jets" in the direction of flow and also indicates the relative weight of forward versus backward transition probability \cite{Meerschaert06,Zhang09x}. The
function $u(x, t)$ can be interpreted as
representing the concentration of a particle plume undergoing anomalous diffusion. The diffusion
coefficient $\xi(x,t)$ depending on both time and space variables satisfies the condition $0 < \xi_{{\rm min}}
\leq \xi(x,t) < \xi_{{\rm max}} < +\infty$, $\forall(x,t)\in[x_L,x_R]
\times[0,t]$, and the forcing function $f(x, t)$ represents the source or sink term. In the current study,
we assume that the problem (\ref{eq1.1}) has a unique and (sufficiently) smooth solution unless otherwise specified \cite{Sandev15,Alikhanov17,Morgado17,Deng2018w}.

The GTSFDE (\ref{eq1.1}) can be regarded as a generalization of classical diffusion equations where
the first-order time derivative is replaced by the generalized Caputo fractional derivative of order
$\gamma \in (0, 1]$ with weighting function $\lambda(t) > 0$ for $t\in [0, T]$, and the second-order spatial derivative is replaced by the
two-sided Riemann-Liouville (R-L) fractional derivative of order $\alpha \in (1, 2]$. Specifically, the
time fractional derivative in Eq. (\ref{eq1.1}) is the generalized Caputo fractional derivative of order
$\gamma$ \cite{Alikhanov17} denoted by
\begin{equation}
{}^{C}_{0}D^{\gamma,\lambda(t)}_{t}u(x,t) = \frac{1}{\Gamma(1 - \gamma)}\int^{t}_0\frac{\lambda(t -
\eta)}{(t - \eta)^{\gamma}}\frac{\partial u(x,\eta)}{\partial \eta}d\eta,
\label{eq1.2}
\end{equation}
which collapses to the widely recognized Caputo or Caputo-tempered fractional derivatives
when $\lambda(t)\equiv1$ or $\lambda(t) = e^{-bt}~(b>0)$~\cite{Xu2016,Chen2018}, respectively. It implies
that the weighting function is indeed often chosen as
$\lambda(t) > 0$ (even with the certain monotonicity) in real-world applications. Meanwhile, the left-handed
(${}_{a}D^{\alpha}_{x}$) and the right-handed (${}_xD^{\alpha}_b$) space fractional derivatives in Eq.
(\ref{eq1.1}) are the R-L fractional derivatives of order $\alpha$ \cite{Podlubny99} which
are defined as
\begin{equation*}
{}_{x_L}D^{\alpha}_{x}u(x,t) = \frac{1}{\Gamma(2 - \alpha)}\frac{\partial^2}
{\partial x^2}\int^{x}_{x_L}\frac{
u(\xi, t)d\xi}{(x - \xi)^{\alpha - 1}}~~\mathrm{and}~~
{}_xD^{\alpha}_{x_R}u(x,t) = \frac{1}{\Gamma(2 - \alpha)}\frac{\partial^2}{\partial x^2}\int^{x_R}_x\frac{u(\xi,
t)d\xi}{(\xi - x)^{\alpha - 1}},
\end{equation*}
where $\Gamma(\cdot)$ denotes the Gamma function. Note that the above equation reduces to the classical
diffusion equation for $\gamma = \lambda(t) \equiv 1$ and $\alpha = 2$.

Generally speaking, although the (semi-)analytical (or closed-form) solutions of particular (generalized) space-time fractional
partial differential equations (PDEs) on the entire real line are accessible via the Laplace or Fourier transforms,
yet these solutions are expressed in terms of special functions which are usually difficult for the numerical
evaluation in practice. Moreover, if we define the problem (\ref{eq1.1}) on a bounded domain, one cannot obtain
any known equations for its fundamental solution; refer to \cite{Mark18,Deng2019}. These naturally promote the
rapid development of numerical methods for fractional PDEs. Therefore, the current study will focus on developing
the numerical approaches for solving the problem (\ref{eq1.1}).

If $\gamma = \lambda(t) \equiv 1$, the problem (\ref{eq1.1}) collapses to the space fractional diffusion
equation (SFDE) with variable coefficients. For such SFDEs, various robust numerical schemes are proposed
by exploiting the shifted Gr\"{u}nwald discretization and the implicit Euler (or Crank-Nicolson)
time-stepping procedure for two-sided R-L fractional derivatives and the first-order time derivative,
respectively; refer to \cite{Meerschaert06,Tadjeran06,Lin2014} for details. To improve the convergence
order of such numerical methods, several studies combined different second-order accurate approximations
for discretizing two-sided R-L fractional derivatives with the Crank-Nicolson technique in order to
obtain the second-order finite difference schemes for solving the SFDEs with variable coefficients. However,
the unconditional convergence of such second-order finite difference schemes is not easy to prove, refer
to \cite{Qu2014x,Sousa15,Feng2015x,Lin2017,Lin2018,Vong19,Lin2019p,Zheng2019,Lin2020} for discussions on this issue. However, these
studies verified the unconditional convergence of second-order finite difference schemes often restrict
diffusion coefficients positively bounded and relied on the spatial variable $x$. Besides, other numerical
treatments including the Chebyshev-tau, finite volume and finite element methods are proposed to solve the
SFDEs with variable coefficients, refer, e.g., to \cite{Ji2012,Ren2013,Doha2013,Ma2014,Liu2014g,Feng2015,Pan2016,Liu2019}
for details.

When $\alpha = 2$, the problem (\ref{eq1.1}) is equivalent to the generalized time fractional diffusion
equation (GTFDE) with variable coefficients. Such GTFDEs were first derived and studied by Sandev \emph{et al.}
in \cite{Sandev15}. Later, Alikhanov adapted the classical \textit{L}1 formula \cite{Podlubny99} and employed
the second-order weighted-shifted Gr\"{u}nwald difference (WSGD) formula \cite{Hao2015} to approximate the generalized Caputo fractional derivative and the spatial R-L fractional
derivative respectively for solving such GTFDEs with variable coefficients. Moreover, the convergence of his implicit
difference schemes is proved to be unconditionally stable, refer to \cite{Alikhanov17} for details. In
addition, Khibie \cite{Khibiev19} has extended Alikhanov's work to establish the stable implicit difference
scheme for solving the multi-term GTFDE with variable coefficients.

On the other hand, although there are several numerical schemes about solving TSFDEs with variable coefficients
--cf. $\lambda(t) \equiv 1$, however those that are proved to be unconditionally convergent  \cite{Firooz,Zhao16,Chi2017,Fu2019,Lin2019n} are only first- and $(2-\gamma)$-order accurate in space and time directions, respectively. It means that proving the unconditional convergence of implicit difference schemes with high-order spatial discretizations is often very challenging. Moreover, there are few results on
numerical solutions of GTSFDEs with variable coefficients via finite difference methods in the literature. Such
GTSFDEs can be regarded as a generalization of the GTFDEs introduced in \cite{Sandev15,Alikhanov17} and their numerical solutions should be more difficult due to lots of computational cost arising from the nonlocal properties in both
spatial and temporal fractional derivatives. Therefore, establishing an unconditionally stable numerical scheme
with low computational cost for solving such GTSFDEs with variable coefficients is a promising topic and also the
main motivation of our current study. In this paper, we develop the implicit difference schemes for GTSFDEs with
variable coefficients, then the implicit schemes are strictly proved to be unconditionally stable and convergent
with second- and $(2-\gamma)$-order accuracy in space and time directions, respectively. Moreover, the implicit
difference schemes lead to the solutions of the resulting linear systems with Toeplitz-like coefficient matrices
which can be solved via direct method in $\mathcal{O}(N^3)$ operations along with $\mathcal{O}(N^2)$ storage.
However, the efficient preconditioned Krylov subspace solvers are employed to reduce the above computational and
memory cost to $\mathcal{O}(N\log N)$ and $\mathcal{O}(N)$, respectively, where $N$ is the number of spatial grid
nodes. Furthermore, the fast sum-of-exponential (SOE) approximation \cite{Jiang17} is extended to reduce computational
and memory cost arising from the nonlocal property in the generalized Caputo fractional derivative with special function
$\lambda(t)$'s. To the best of our knowledge, this is the first successful attempt to derive such a fast and stable
numerical scheme of GTSFDEs with variable coefficients. Meanwhile, numerical experiments are reported to support
our theoretical finding and effectiveness of the proposed schemes.

The rest of this paper is organized as follows. In Section \ref{sec2}, the approximations of the generalized
Caputo and R-L fractional derivatives are recalled to establish the implicit difference scheme. Meanwhile, the
stability and convergence of the proposed difference scheme are proved in details. In Section \ref{sec3}, the
practical implementation of the proposed schemes requires to solve a sequence of linear systems with Toeplitz-like coefficient matrices. The efficient preconditioned Krylov subspace solvers are adapted and investigated to
handle such Toeplitz-like resultant linear systems. In Section \ref{sec4}, numerical experiments are reported
to demonstrate the efficiency of the proposed method. Some concluding remarks are given in Section \ref{sec5}.
\section{An implicit difference scheme for GTSFDEs}
\label{sec2}
In this section, we first review the approximation of the generalized Caputo fractional
derivative and employ the second-order WSGD approximation \cite{Hao2015} to derive
the implicit difference scheme to problem (\ref{eq1.1}). Moreover, we
have to provide the certain smoothness and monotonicity \cite{Alikhanov17} for the weighting function
$\lambda(t)$, then we can derive in details both the stability and convergence of our implicit difference scheme.
\subsection{The approximation for the generalized Caputo fractional derivative}
\label{sec2.1}
We first briefly recall the generalized $L1$ formula for approximating the temporal
fractional derivative ${}^{C}_{0}D^{\gamma,\lambda(t)}_{t}$ proposed in \cite{Alikhanov17}
and denote its approximation result by $\Delta^{\gamma,\lambda(t)}_{0,t}$. To derive
the difference scheme, we first introduce a rectangle $\bar{Q}_T = \{(x,t):~x_L
\leq x \leq x_R,~0\leq t \leq T\}$ discretized on the mesh $\varpi_{h,\tau}
= \varpi_h\times\varpi_{\tau}$, where $\varpi_h = \{x_i = x_L + ih,~0\leq i
\leq N,~h = \frac{x_R - x_L}{N}\}$ and $\varpi_{\tau}
= \{t_j = j\tau,~j=0,1,\ldots,M,~\tau = \frac{T}{M}\}$. We also denote by ${\bm v} = \{v_i~|~i
= 0,1,\ldots,N\}$ any grid function.
%
%
Moreover, we denote the linear interpolation over
the time interval $(t_j,t_{j+1})$ with $0\leq j\leq M - 1$ by
\begin{equation*}
\Pi_{1,s}v(t) = v(t_{s + 1})\frac{t - t_s}{\tau} + v(t_s)\frac{t_{s+1} - t}{\tau}.
\end{equation*}

At each time step $t_{j + 1}$ with $j=0,1,\ldots,M-1$, the generalized $L1$ formula is defined by
\begin{equation*}
\begin{split}
{}^{C}_{0}D^{\gamma,\lambda(t)}_{t}v(t)|_{t = t_{j + 1}} & = \frac{1}{\Gamma(1 - \gamma)}
\int^{t_{j+1}}_{0}\frac{\lambda(t_{j+1} - \eta)v'(\eta)d\eta}{(t_{j+1} - \eta)^{\gamma}}\\
& = \frac{1}{\Gamma(1 - \gamma)}\left[\sum^{j}_{s = 0}v_{t,s}\int\limits^{t_{s+1}}_{t_s}\frac{
\lambda(t_{j+1} - \eta)d\eta}{(t_{j+1} - \eta)^{\gamma}} + \sum^{j}_{s = 0}\int\limits^{t_{s+1}
}_{t_s}\frac{\lambda(t_{j + 1} - \eta)[v(\eta) - \Pi_{1,s}v(\eta)]'d\eta}{(t_{j + 1} -
\eta)^{\gamma}}\right]\\
& = \frac{\tau^{1 - \gamma}}{\Gamma(2 - \gamma)}\sum\limits^{j}_{s = 0}[\lambda_{j-s + 1/2}
a_{j - s} + (\lambda_{j - s} - \lambda_{j - s + 1})b_{j - s}]v_{t,s} + R^{j}_1 + R^{j}_2,
\end{split}
\end{equation*}
where $\lambda_s = \lambda(t_s)$ and
\begin{equation*}
v_{t,s} = \frac{v(t_{s+1})-v(t_s)}{\tau},~~a_{\ell} = (\ell + 1)^{1 - \gamma} - \ell^{1 - \gamma},~~b_{\ell} =
\frac{1}{2 - \gamma}[(\ell + 1)^{2 - \gamma} - \ell^{2 - \gamma}] - \frac{1}{2}[(\ell + 1)^{1
-\gamma} + \ell^{1 - \gamma}],\quad \ell\geq 1,
\end{equation*}
and the definition of $R^{j}_1,~R^{j}_2$ and their estimations can be separately found in
\cite{Alikhanov17}. 
The truncation error and property of the generalized $L1$ formula are also analyzed in
\cite[Lemma 4.1]{Alikhanov17} as follows
%
\begin{lemma} Assume that $\gamma\in(0,1)$, $\lambda(t)>0$, $\lambda'(t)\leq 0$, and $\lambda(t),
v(t)\in\mathcal{C}^2[0,t_{j+1}]$. Then
\begin{equation}
{}^{C}_{0}D^{\gamma,\lambda(t)}_{t}v(t_{j+1}) = \Delta^{\gamma,\lambda(t)}_{0,t_{j+1}}v^{j+1} +
\mathcal{O}(\tau^{2 - \gamma}),
\label{eq1.3}
\vspace{-6mm}
\end{equation}
where $\Delta^{\gamma,\lambda(t)}_{0,t_{j+1}}u^{j+1} = \sum\limits^{j}_{s = 0}c_{j-s}[u(t_{s + 1}) -
u(t_{s})]$ and $c_{\ell} = \frac{\tau^{-\gamma}}{\Gamma(2 - \gamma)}[\lambda_{\ell + 1/2}a_{\ell
} + (\lambda_{\ell} - \lambda_{\ell + 1})b_{\ell}]~(\ell \geq 0)$. Moreover, the following inequalities hold:
\begin{equation*}
a_0 > a_1 > \cdots > a_{\ell} > \frac{1 - \gamma}{(\ell + 1)^{\gamma}},\quad b_0 > b_1 > \cdots > b_{\ell} > 0.
\end{equation*}
\label{lem2.1}
\vspace{-3mm}
\end{lemma}
\vspace{-3mm}
Based on the property of $a_{\ell}$ and $b_{\ell}$, we can obtain the following result for the
coefficients $c_{\ell}$, which is absolutely vital for our theoretical analysis in the next subsection.
\begin{lemma} For all $\ell = 0, 1,\ldots,~\gamma\in(0, 1)$ and $\lambda(t)\in \mathcal{C}^2[0, T]$,
where $\lambda(t) > 0,~\lambda'(t) \leq 0$ for all $t\in[0, T]$, the following inequalities hold:
\begin{equation*}
c_0 > c_1 > \ldots > c_{\ell} > \frac{\lambda(t_{\ell + 1/2})}{\Gamma(1 - \gamma)t^{\gamma}_{\ell + 1}}.
\end{equation*}
\vspace{-3mm}
\end{lemma}
After we introduce the temporal discretization, it is the time to characterize the discretization in the
space variable. First of all, we denote by
\begin{equation*}
\mathcal{L}^{n + \alpha}(\mathbb{R}) = \left\{v \mid v\in L_{1}(\mathbb{R})~{\rm and}
~\int^{+\infty}_{-\infty}(1 + |k|)^{n +\linebreak {\color{blue}\alpha}}|\hat{v}(k)|dk < \infty\right\},
\end{equation*}
where $\hat{v}(k) = \int^{+\infty}_{-\infty}e^{\iota kx}v(x)dx$ is
the Fourier transformation of $v(x)$, and by $\iota = \sqrt{-1}$ the imaginary unit. Then we introduce
the following preliminary lemma, which provides numerical approximations for the spatial R-L fractional
derivatives:
\begin{lemma} Let $v(x)\in\mathcal{L}^{2 + \alpha}(\mathbb{R})$ and define the following
difference operators
\begin{equation*}
\delta^{\alpha}_{x,+}v(x) = \frac{1}{h^{\alpha}}\sum^{[\frac{x - x_L}{h}]}_{k = 0}w^{(\alpha)}_kv(x - (k
- 1)h)~~{\rm and}~~
\delta^{\alpha}_{x,-}v(x) = \frac{1}{h^{\alpha}}\sum^{[\frac{x_R - x}{h}]}_{k = 0}w^{(\alpha)}_kv(x + (k
- 1)h).
\end{equation*}
Then, for a fixed $h$, we have
\begin{equation*}
{}_aD^{\alpha}_xv(x) = \delta^{\alpha}_{x,+}v(x) + \mathcal{O}(h^2)\quad\quad{\rm and}\quad\quad {}_x
D^{\alpha}_xv(x) = \delta^{\alpha}_{x,-}v(x) + \mathcal{O}(h^2),
\end{equation*}
where $[\cdot]$ is the floor function and
\begin{equation*}
\begin{cases}
w^{(\alpha)}_0 = \kappa_1g^{(\alpha)}_0,\quad w^{(\alpha)}_1 = \kappa_1g^{(\alpha)}_1 + \kappa_0g^{(0)
}_0,\\
w^{(\alpha)}_k = \kappa_1g^{(\alpha)}_k + \kappa_0g^{(\alpha)}_{k-1} + \kappa_{-1}g^{(\alpha)}_{k-2}, &
k \geq 2,
\end{cases}
\end{equation*}
with
\begin{equation*}
\kappa_1 = \frac{\alpha^2 + 3\alpha + 2}{12},~~\kappa_0 = \frac{4 - \alpha}{6},~~\kappa_{-1} = \frac{
\alpha^2 - 3\alpha + 2}{12},~~{\rm and}~~g^{(\alpha)}_k = (-1)^k\binom{\alpha}{k}.
\end{equation*}
\label{lem2.2}
\vspace{-2mm}
\end{lemma}
\vspace{-2mm}
At this stage, the numerical approximations of both the temporal and spatial
fractional derivatives have been set for the derivation of the targeted implicit difference
scheme. Let $u(x, t)\in \mathcal{C}^{4,2}_{x,t}([x_L, x_R]\times[0, T])$ be a solution to the problem
(\ref{eq1.1}). Then we consider Eq. (\ref{eq1.1}) at the set of grid points $(x, t) = (x_i, t_{j+1})\in
\bar{Q}_{T},~i = 1, 2,\ldots,N-1,~j = 0, 1,\ldots, M-1$:
\begin{equation*}
{}^{C}_0D^{\gamma,\lambda(t)}_{t}u(x_i,t_{j+1}) = \xi(x_i,t_{j+1})\Big[p{}_{x_L}D^{\alpha}_xu(x,t)
+ (1 - p){}_{x}D^{\alpha}_{x_R}u(x,t)\Big]_{(x_i,t_{j + 1})} + f(x_i,t_{j+1}).
\end{equation*}
Let $U$ be a grid function defined by
\begin{equation*}
U^{j}_i: = u(x_i,t_j)\quad{\rm and}\quad f^{j}_i = f(x_i,t_j),\quad 0\leq i\leq N,\quad 0\leq j\leq M.
\end{equation*}
Using this notation and recalling Lemma \ref{lem2.1} and Lemma \ref{lem2.2}, we can write the problem (\ref{eq1.1})
at the grid points $(x_i, t_{j+1})$ as follows
\begin{equation}
\Delta^{\gamma,\lambda(t)}_{0,t_{j+1}}U^{j+1}_{i} = \xi^{j+1}_{i}(\delta^{\alpha}_h U^{j+1}_i)
+ f^{j+1}_i + R^{j + 1}_i,\quad 1\leq i \leq N-1,\quad 0\leq j\leq M-1,
\end{equation}
where $\{R^{j+1}_i\}$ are small and satisfy the relation $|R^{j+1}_i| = \mathcal{O}(\tau^{2 - \gamma} + h^2)$
for $1\leq i\leq N - 1,~0\leq j\leq M-1$. We omit them and use the initial-boundary value conditions
\begin{equation*}
\begin{cases}
U^{0}_i = \phi(x_i), & 1\leq i\leq N-1,\\
U^{j}_0 = \varphi(t_{j}),\quad U^{j}_N = \psi(t_{j}), & 0\leq j\leq M.
\end{cases}
\end{equation*}
%
%

%
For the sake of clarity, we introduce the notations
\begin{equation*}
\xi^{j}_{i} = \xi(x_{i},t_{j}),\quad
\delta^{\alpha}_hu^{j+1}_{i}  = \frac{1}{h^{\alpha}}\left[p\sum^{i + 1}_{k = 0}w^{(\alpha)}_ku^{j+1}_{i -
k + 1} + (1-p)\sum^{N - i + 1}_{k = 0}w^{(\alpha)}_ku^{j+1}_{i + k - 1}\right],
\end{equation*}
and then we arrive at the implicit difference scheme with (local) truncation errors of
$\mathcal{O}(\tau^{2 - \gamma} + h^2)$:
\begin{equation}
\begin{cases}
\Delta^{\gamma,\lambda(t)}_{0,t_{j+1}}u^{j+1}_{i} = \xi^{j+1}_{i}(\delta^{\alpha}_h u^{j+1}_i)
+ f^{j+1}_i,& i = 1,2,\ldots,N-1,\quad j = 0,1,\ldots, M-1,\\
u^{0}_i = \phi(x_i),& i = 0,1,\ldots,N,\\
u^{j}_0 = \varphi(t_{j}),\quad u^{j}_N = \psi(t_{j}),& j = 0,1,\ldots,M.
\end{cases}
\label{eq2.3}
\end{equation}
It is interesting to note that for $\lambda(t) \equiv 1$ and $\gamma\rightarrow 1$, Eq. (\ref{eq2.3})
reduces to the classical backward Euler scheme for solving the SFDEs with variable coefficients \cite{Meerschaert06,Lei2013x}.
Similarly, if $\alpha = 2$, the above scheme (\ref{eq2.3}) collapses to the implicit difference scheme
introduced in \cite{Alikhanov17} for solving the variable-coefficient GTFDEs.
\subsection{Stability and convergence analysis}
\label{sec2.2}
In this subsection, we are committed to analyzing both the stability and convergence for the implicit difference
scheme (\ref{eq2.3}). We define
\begin{equation*}
V_h = \{{\bm v}~|~{\bm v} = \{v_i\}{\rm~is~a~grid~function~on}~\varpi_h~{\rm and}~v_i = 0~{\rm if}~i = 0, N\},
\end{equation*}
and, for all ${\bm u}, {\bm v}\in V_h$, the discrete inner product and corresponding discrete $L^2$-norms
\begin{equation*}
({\bm u}, {\bm v}) = h\sum^{N-1}_{i=1}u_iv_i,~~~and~~~\|{\bm u}\| = \sqrt{({\bm u}, {\bm u})}.
\end{equation*}
\quad~The starting point of our analysis is the following theoretical result.
\begin{lemma} {\rm (\cite{Vong2016,Gu2017})}
Let $\alpha\in(1,2)$ and $g^{(\alpha)}_k$ be defined in Lemma \ref{lem2.2}, then we obtain
\begin{equation*}
\begin{cases}
w^{(\alpha)}_0 = \kappa_1 > 0,\quad w^{(\alpha)}_1 < 0,\quad w^{(\alpha)}_k > 0,\quad k\geq 3, \\
\sum\limits^{\infty}_{k = 0} w^{(\alpha)}_k = 0,\quad \sum\limits^{N}_{k = 0}w^{(\alpha)}_k
< 0,\quad N > 1,\\
w^{(\alpha)}_0 + w^{(\alpha)}_2 \geq 0.
\end{cases}
\vspace{-5mm}
\end{equation*}
\label{lem2.3}
\vspace{-5mm}
\end{lemma}

\begin{figure}[!htpb]
\centering
\includegraphics[width=3.06in,height=2.75in]{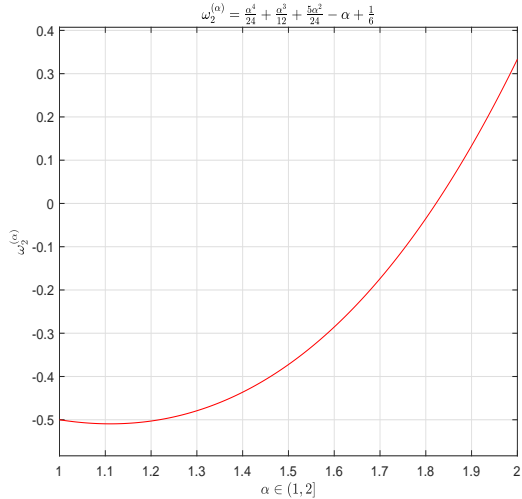}
\caption{The plot of $w^{(\alpha)}_2$ with $\alpha\in(1,2]$.}
\label{fig1}
\end{figure}

In fact, this lemma does not show whether $w^{(\alpha)}_2$ is positive or negative. After simple calculations,
we obtain
\begin{equation}
\begin{split}
w^{(\alpha)}_2 &= \kappa_1g^{(\alpha)}_2 + \kappa_0g^{(\alpha)}_1 + \kappa_{-1}g^{(\alpha)}_0\\
& = \frac{\alpha^4}{24} + \frac{\alpha^3}{12} + \frac{5\alpha^2}{24} - \alpha + \frac{1}{6},
\end{split}
\end{equation}
where $\alpha\in(1,2]$ and it can be plotted as in Fig. \ref{fig1}.
As seen from Fig. \ref{fig1}, the following proposition can be derived, which is helpful to analyse
the property of the coefficient matrices appearing in Eq. (\ref{eq1.9a}) in the next section.
\begin{proposition}
When $\alpha\in(1,\alpha_0)$, then $w^{(\alpha)}_2 < 0$. Similarly, when $\alpha\in[\alpha_0,2]$,
then $w^{(\alpha)}_2 \geq 0$ with $\alpha_0 \approx 1.8223$. Moreover, the sufficient condition for
$W_{\alpha}$ and $W^{T}_{\alpha}$ to be diagonally dominant is $\alpha\in[\alpha_0,2]$,
where the matrix
\begin{equation}
W_{\alpha} = \begin{pmatrix}
w^{(\alpha)}_1     &w^{(\alpha)}_0     & 0               & \cdots   & 0 &0  \\
w^{(\alpha)}_2     &w^{(\alpha)}_1     & w^{(\alpha)}_0  & 0 & \cdots &0 \\
\vdots             &w^{(\alpha)}_2     & w^{(\alpha)}_1  & \ddots & \ddots & \vdots \\
\vdots             &\ddots             & \ddots          &    \ddots & \ddots  &  0                 \\
w^{(\alpha)}_{N-2} &\ddots             & \ddots          & \ddots & w^{(\alpha)}_1 & w^{(\alpha)}_0\\
w^{(\alpha)}_{N-1} &w^{(\alpha)}_{N-2} & \cdots          & \cdots & w^{(\alpha)}_2 & w^{(\alpha)}_1
\end{pmatrix}\in \mathbb{R}^{(N - 1)\times (N-1)}.
\label{eq2.4add}
\end{equation}
\label{rem2.1}
\vspace{-2mm}
\end{proposition}
\noindent\textbf{Proof}. Since $\alpha\in[\alpha_0,2]$, it holds $w^{(\alpha)}_1 < 0$ and $w^{(\alpha)}_k
\geq 0$ ($k\neq 1$). According to $\sum\limits^{\infty}_{k = 0}w^{(\alpha)}_k = 0$, it holds that both $W_{\alpha}$
and $W^{T}_{\alpha}$ are diagonally dominant \cite{Lin2020}. \hfill$\Box$

Based on Lemma \ref{lem2.3}, the first two properties of the discrete inner product related to
two approximate operators $\delta^{\alpha}_{x,+}$ and $\delta^{\alpha}_{x,-}$ can be shown
below.
\begin{lemma}{\rm(\cite{Vong2016,Gu2017})}
For $\alpha \in (1,2)$ and $N \geq 5$, and any ${\bm v}\in V_h$, it holds that
\begin{equation*}
(-\delta^{\alpha}_{x,+}{\bm v},{\bm v}) = (-\delta^{\alpha}_{x,-}{\bm v},{\bm v}) > c\ln 2\|{\bm v}\|^2.
\end{equation*}
where $c$ is positive constant independent of the spatial step size $h$.
\label{lem2.5}
\end{lemma}
\begin{theorem}
For $\alpha \in (1, 2)$, and any ${\bm v}\in V_h$, it holds that
\begin{equation*}
(\delta^{\alpha}_h{\bm v},{\bm v}) < -c\ln 2\|{\bm v}\|^2,
\end{equation*}
where $c$ is the same constant appearing in Lemma \ref{lem2.5}.
\label{lem2.6}
\end{theorem}
\textbf{Proof}. The concrete expression of $(\delta^{\alpha}_h{\bm v},{\bm v})$ can
be written as
\begin{equation*}
(\delta^{\alpha}_h{\bm v},{\bm v}) = p(\delta^{\alpha}_{x,+}{\bm v},{\bm v})
+ (1-p)(\delta^{\alpha}_{x,-}{\bm v},{\bm v})\leq -c\ln{2}\|{\bm v}\|_{2},
\end{equation*}
and this completes the proof of Theorem \ref{lem2.6}.\hfill $\Box$

To establish the stability of the difference scheme, we still need to introduce the following lemma.
\begin{lemma}
For any function $v(t)$ defined on the discrete grid $\varpi_{\tau}=\{t_j = j\tau:j = 0, 1,\ldots, M\}$,
the following inequality holds
\begin{equation}
{\bm v}^{j+1}(K^{j+1})^{-1}\Delta^{\gamma,\lambda(t)}_{0,t_{j+1}}{\bm v} \geq \frac{1}{2}\Delta^{\gamma,
\lambda(t)}_{0,t_{j+1}}\|{\bm v}\|^{2}_{(K^{j+1})^{-1}},
\end{equation}
where $K^{j+1} = {\rm diag}(\xi^{j+1}_1,\xi^{j+1}_2,\ldots,\xi^{j+1}_{n-1}) > 0$ and
$\|{\bm v}\|^{2}_{(K^{j+1})^{-1}} = {\bm v}^T(K^{j+1})^{-1}{\bm v}$.
\label{lem2.5x}
\end{lemma}
\noindent\textbf{Proof}. We rewrite the following inner product
\begin{equation}
{\bm v}^{j+1}(K^{j+1})^{-1}\Delta^{\gamma,\lambda(t)}_{0,t_{j+1}}{\bm v} = \tilde{{\bm v}}^{j+1}\Delta^{\gamma,
\lambda(t)}_{0,t_{j+1}}\tilde{{\bm v}} \geq \frac{1}{2}\Delta^{\gamma,\lambda(t)}_{0,t_{j+1}}\|\tilde{{\bm v}}\|_2,
\label{ineq2.6}
\end{equation}
where $\tilde{{\bm v}} = (K^{j+1})^{-\frac{1}{2}}{\bm v}$ regarded as a (weighted) function $v(t)$ defined
on the discrete grid $\varpi_{\tau}$. Meanwhile, the inequality (\ref{ineq2.6}) is correct due to
\cite[Lemma 4.4]{Alikhanov17}. \hfill $\Box$

Another ingredient, introduced as the following lemma, is also required to describe the diagonally weighted
norm that will be used in the next theorem.
\begin{lemma}{\rm(\cite{Lin2019p})}
Let $H\in\mathbb{R}^{n\times n}$ be a symmetric matrix with eigenvalues $\tilde{\lambda}_1
\geq \tilde{\lambda}_2 \geq \ldots \geq \tilde{\lambda}_n$. Then for all ${\bm w}\in\mathbb{R}^{n\times 1}$,
\begin{equation}
\tilde{\lambda}_n{\bm w}^T{\bm w}\leq {\bm w}^T H{\bm w} \leq \tilde{\lambda}_1{\bm w}^T{\bm w}.
\end{equation}
\vspace{-3mm}
\label{lem2.8}
\end{lemma}
\vspace{-4mm}
Now we can conclude the stability and convergence of the implicit difference
scheme (\ref{eq2.3}). For simplicity of presentation, we denote $a^{j+1}_s =
c_{j-s}$, then $\Delta^{\gamma,\lambda(t)}_{0,t_{j+1}}{\bm u} = \sum\limits^{
j}_{s = 0} ({\bm u}^{s+1} - {\bm u}^{s})a^{j+1}_{s}$.
\vspace{-2mm}
\begin{theorem}
If we define $\|{\bm f}^{j+1}\|^2 = h\sum\limits^{N-1}_{i=1}f^2(x_i, t_{j+1})$, then the implicit difference
scheme (\ref{eq2.3}) is unconditionally stable and the following a priori estimate holds:
\begin{equation}
\|{\bm u}^{j+1}\|^{2}_{(K^{j+1})^{-1}} \leq \frac{1}{\xi_{\rm min}}\left[\|{\bm u}^0\|^{2} +
\frac{\Gamma(1 - \gamma)T^{\gamma}}{2c\xi_{{\rm min}}\ln2\lambda(T)}\max\limits_{0\leq j\leq M-1}\|{\bm f}^{j+1}\|^{2}\right]
\label{eq2.5}
\end{equation}
where ${\bm u}^{j+1} = [u^{j+1}_1, u^{j+1}_2,\ldots, u^{j+1}_{N-1}]^T$.
\label{thm2.2}
\end{theorem}
\textbf{Proof}. To make an inner product of Eq. (\ref{eq2.3}) with ${\bm u}^{j+1}$, we have
\begin{equation}
(\Delta^{\gamma,\lambda(t)}_{0,t_{j+1}}{\bm u},(K^{j+1})^{-1}{\bm u}^{j+1}) = (\delta^{\alpha}_h {\bm u}^{j+1},{\bm u}^{j+1})
+ ({\bm f}^{j+1},(K^{j+1})^{-1}{\bm u}^{j+1}).
\label{eq2.4x}
\end{equation}
It follows from Theorem \ref{lem2.6} and Lemma \ref{lem2.5x} that
\begin{equation}
(\delta^{\alpha}_h {\bm u}^{j+1},{\bm u}^{j+1}) \leq -c\ln2\|{\bm u}^{j+1}\|^2
\label{eq2.5x}
\end{equation}
and
\begin{equation}
(\Delta^{\gamma,\lambda(t)}_{0,t_{j+1}}{\bm u},(K^{j+1})^{-1}{\bm u}^{j+1}) \geq \frac{1}{2}
\Delta^{\gamma,\lambda(t)}_{0,t_{j+1}}\|{\bm u}\|^{2}_{(K^{j+1})^{-1}}.
\label{eq2.6x}
\end{equation}
Substituting (\ref{eq2.5x})-(\ref{eq2.6x}) into (\ref{eq2.4x}) and using
the Cauchy-Schwarz and Young's inequalities, we obtain
\begin{equation*}
\begin{split}
\frac{1}{2}
\Delta^{\gamma,\lambda(t)}_{0,t_{j+1}}\|{\bm u}\|^{2}_{(K^{j+1})^{-1}} &\leq -c\ln2\|{\bm u}^{j+1}\|^2
+ ({\bm f}^{j+1},(K^{j+1})^{-1}{\bm u}^{j+1})\\
&\leq -c\ln2\|{\bm u}^{j+1}\|^2 + c\xi_{{\rm min}}\ln2\|{\bm u}^{j+1}\|^{2}_{(K^{j+1})^{-1}} + \frac{1}{4c
\xi_{{\rm min}}\ln2}\|{\bm f}^{j+1}\|^{2}_{(K^{j+1})^{-1}}\\
&\leq -c\ln2\|{\bm u}^{j+1}\|^2 + c\ln2\|{\bm u}^{j+1}\|^{2} + \frac{1}{4c\xi_{{\rm min}}\ln2}\|{\bm f}^{j+1}
\|^{2}_{(K^{j+1})^{-1}}\quad {\rm(cf.~Lemma \ref{lem2.8})} \\
& = \frac{1}{4c\xi_{{\rm min}}\ln2}\|{\bm f}^{j+1}\|^{2}_{(K^{j+1})^{-1}}.
\end{split}
\end{equation*}
\quad~Next, we have the following inequality
\begin{equation}
a^{j+1}_j\|{\bm u}^{j+1}\|^{2}_{(K^{j+1})^{-1}} \leq \sum^{j}_{s = 1}(a^{j+1}_s - a^{j+1}_{s-1})\|{\bm u}^s
\|^{2}_{(K^{j+1})^{-1}} + a^{j+1}_0\|{\bm u}^0\|^{2}_{(K^{j+1})^{-1}} + \frac{1}{2c\xi_{{\rm min}}\ln2}\|{\bm
f}^{j+1}\|^{2}_{(K^{j+1})^{-1}}.
\end{equation}
Employing the inequality $a^{j+1}_0 = c_j > \frac{\lambda(T)}{\Gamma(1 -
\gamma)T^{\gamma}}$ (cf. \cite[Theorem 5.1]{Alikhanov17}), we obtain
\begin{equation}
\begin{split}
a^{j+1}_j\|{\bm u}^{j+1}\|^{2}_{(K^{j+1})^{-1}} & \leq \sum^{j}_{s = 1}(a^{j+1}_s - a^{j+1}_{s-1})\|{\bm u}^s
\|^{2}_{(K^{j+1})^{-1}}\\
&\quad + a^{j+1}_0\left[\|{\bm u}^0\|^{2}_{(K^{j+1})^{-1}} + \frac{\Gamma(1 - \gamma)T^{\gamma}}{2c\xi_{{\rm min}}
\ln2\lambda(T)}\|{\bm f}^{j+1}\|^{2}_{(K^{j+1})^{-1}}\right].
\end{split}
\label{eq2.11}
\end{equation}
\quad~Suppose $h < 1$ and denote
\begin{equation*}
\mathcal{P} \triangleq \frac{1}{\xi_{\rm min}}\left[\|{\bm u}^0\|^{2} + \frac{\Gamma(1 - \gamma)T^{\gamma}}{2c
\xi_{{\rm min}}\ln2\lambda(T)}\max\limits_{0\leq j\leq M-1}\|{\bm f}^{j+1}\|^{2}\right].
\end{equation*}
Then, Eq. (\ref{eq2.11}) can be rewritten as
\begin{equation}
a^{j+1}_j\|{\bm u}^{j+1}\|^{2}_{(K^{j+1})^{-1}} \leq \sum^{j}_{s = 1}(a^{j+1}_s - a^{j+1}_{s-1})\|{\bm u}^s\|^{2}_{(K^{j+1})^{-1}}
+ a^{j+1}_0\mathcal{P}.
\end{equation}
\quad~At this stage, by mathematical induction we prove that
\begin{equation}
\|{\bm u}^s\|^{2}_{(K^{j+1})^{-1}}\leq \mathcal{P},\quad~0\leq s\leq j+1,
\label{eq2.13}
\end{equation}
is valid for the fixed $j$. The result is obviously true for $s =
0$ from (\ref{eq2.11}). Assuming that (\ref{eq2.13}) holds for all $0 \leq s \leq j~(0 \leq j \leq
M - 1)$, then from (\ref{eq2.11}) at $0 \leq s \leq j+1$, one has
\begin{equation*}
\begin{split}
a^{j+1}_j\|{\bm u}^{j+1}\|^{2}_{(K^{j+1})^{-1}} & \leq \sum^{j}_{s = 1}(a^{j+1}_s - a^{j+1}_{s-1})\|{\bm u}^s\|^{2}_{(K^{j+1})^{-1}}
+ a^{j+1}_0\mathcal{P}\\
& \leq \sum^{j}_{s = 1}(a^{j+1}_s - a^{j+1}_{s-1})\mathcal{P} + a^{j+1}_0\mathcal{P}\\
& = a^{j+1}_j\mathcal{P},
\end{split}
\end{equation*}
This completes the proof of Theorem \ref{thm2.2}. \hfill $\Box$

The following theorem shows that our proposed implicit difference scheme achieves $(2-\gamma)$-order
and quadratic-order convergence in time and space variables, respectively, when the solution of Eq.
(\ref{eq1.1}) is sufficiently smooth. To our knowledge, it is the first theoretical result on the
convergence of implicit difference schemes for solving the variable-coefficient GTSFDEs (\ref{eq1.1}).
\begin{theorem}
Suppose that $u(x, t)\in \mathcal{C}^{4,2}_{x,t}([x_L, x_R]\times[0, T])$ is the solution of
Eq. (\ref{eq1.1}) and $\{u^{j}_i|x_i\in \varpi_h,~ 0 \leq j \leq M\}$ is the
solution of the implicit difference scheme (\ref{eq2.3}). Define
\begin{equation}
E^{j}_i = u(x_i,t_j) - u^{j}_i,\quad x_i\in \varpi_h,~ 0 \leq j \leq M,
\end{equation}
where $\varpi_h = \{x_i = i h, i =0, 1,\ldots,N;~Nh = b - a\}$, then there exists
a positive constant $\tilde{c}$ such that
\begin{equation*}
\|E^{j}\| \leq \tilde{c} (\tau^{2 - \gamma} + h^2),\quad 0 \leq j \leq M.
\end{equation*}
\label{thm2.3}
\vspace{-8mm}
\end{theorem}
\textbf{Proof}. It can be easily obtained that $E^j$ satisfies the following
error equation
\begin{equation}
\begin{cases}
\Delta^{\gamma,\lambda(t)}_{0,t_{j+1}}E^{j+1}_{i} = \xi^{j+1}_{i}(\delta^{
\alpha}_h E^{j+1}_i) + R^{j+1}_i,& i = 1,2,\ldots,N-1,
\quad j = 0,1,\ldots, M-1,\\
E^{0}_i = 0,& i = 0,1,\ldots,N,\\
E^{j}_0 = 0,\quad E^{j}_N = 0,& j = 0,1,\ldots,M,
\end{cases}
\end{equation}
where ${\bm R}^{j+1} = [R^{j+1}_1,R^{j+1}_2,\cdots,R^{j+1}_{N-1}]^T$ and the truncation error term is $\|{\bm R}^{j+1}\|
= \mathcal{O}(\tau^{2 - \gamma} + h^2)$. In virtue of Theorem
\ref{thm2.2} and Lemma \ref{lem2.8}, we define ${\bm E}^{j+1} = [E^{j+1}_1,E^{j+1}_2,\cdots,E^{j+1}_{N-1}]^T$ and
then arrive at
\begin{equation*}
\|{\bm E}^{j+1}\|^{2}_{(K^{j+1})^{-1}}\leq \frac{\Gamma(1 - \gamma)T^{\gamma}}{2c\xi_{{\rm min}}\ln2\lambda(T)}\|
{\bm R}^{j+1}\|^{2}_{(K^{j+1})^{-1}}
\Rightarrow \|{\bm E}^{j+1}\| \leq \tilde{c}(\tau^{2 - \gamma} + h^2),\quad 0\leq j\leq M-1,
\end{equation*}
which proves the theorem.\hfill $\Box$

Theorem \ref{thm2.3} implies that our numerical scheme converges to the optimal order
$\mathcal{O}(\tau^{2 - \gamma} + h^2)$ in the $L^2$-norm, when the solution of Eq.
(\ref{eq1.1}) is sufficiently smooth. Besides, if the solution of Eq. (\ref{eq1.1}) is
non-smooth, several useful alternatives utilizing the non-uniform temporal step or initial
correction techniques \cite{Jin2016x,Jin2019x,Stynes17,Liao18} can be adapted to
address this problem. However, that is not the emphasis of this current study
and we point the reader to the next section for a short
discussion. In addition, the above analysis can be
similarly adapted to remedy defects in our previous work \cite{Gu2017}, which only focuses on
the model problem with time-varying diffusion coefficients.
\section{Efficient implementation of the proposed implicit difference scheme}
\label{sec3}
In order to develop an efficient implementation of the proposed scheme, we rewrite the implicit difference
scheme (\ref{eq2.3}) into the following form with $i = 1,2,\ldots,N-1$
and $j = 0,1,\ldots,M-1$:
\begin{equation}
(c_0 u^{j+1}_{i} - c_ju^{0}_i) - \sum^{j}_{s = 1}(c_{s - 1} - c_s)u^{j + 1 - s}_{i}
= \frac{\xi^{j+1}_{i}}{h^{\alpha}}\left[p\sum^{i + 1}_{k = 0}w^{(\alpha)}_{k}u^{j + 1}_{i - k + 1} + (1 -
p)\sum^{N - i + 1}_{k = 0}w^{(\alpha)}_k u^{j+ 1}_{i + k - 1}\right] + f^{j + 1}_i,
\end{equation}
or, equivalently,
\begin{equation}
c_0 u^{j+1}_{i} - \frac{\xi^{j+1}_{i}}{h^{\alpha}}\left[p\sum^{i + 1}_{k = 0}w^{(\alpha)}_{k}u^{j +
1}_{i - k + 1} + (1 - p)\sum^{N - i + 1}_{k = 0}w^{(\alpha)}_k u^{j+ 1}_{i + k - 1}\right] = c_ju^{0}_i
+ \sum^{j}_{s = 1}(c_{s - 1} - c_s)u^{j + 1 - s}_{i} + f^{j + 1}_i.
\label{eq3.1}
\end{equation}

At this stage, the above implicit difference scheme can be reformulated as the following sequence of linear systems,
\begin{equation}
\mathcal{M}^{(j + 1)}{\bm u}^{j + 1} = c_j{\bm u}^{0} + \sum\limits^{j}_{s = 1}(c_{s - 1} - c_{s}){\bm u}^{j + 1 - s}
+ {\bm f}^{j + 1},\quad j = 0,1,2,\ldots,M-1,
\label{eq1.9a}
\end{equation}
where $\mathcal{M}^{j+1} = c_0I - \frac{K^{(j+1)}}{h^{\alpha}}\left[pW_{\alpha} + (1 - p)W^{T}_{\alpha}\right]$, ${\bm u
}^{j} = [u^{j}_1,u^{j}_2,\ldots,u^{j}_{N-1}]^T$, ${\bm f}^{j} = [f^{j}_1,f^{j}_2,\ldots,f^{j}_{N-1}]^T$, $K^{(j+1)} =
{\rm diag}(\xi^{j+1}_1,\xi^{j+1}_2,\cdots,\xi^{j+1}_{N-1})$ and $I$ is the identity matrix of order $(N-1)$. Meanwhile,
it is obvious that $W_{\alpha}$ (\ref{eq2.4add}) is a Toeplitz matrix; refer to \cite{Meerschaert06,Ngbook04}. Therefore, it can be stored
with $N$ entries and the matrix-vector product involving the matrix $\mathcal{M}^{(j)}$ can be evaluated via fast Fourier
transforms (FFTs) in $\mathcal{O}(N\log N)$ operations \cite{Ngbook04,Lei2013x}. On the other hand, it is meaningful to
remark that the sequence of linear systems (\ref{eq1.9a}) corresponding to the implicit scheme (\ref{eq2.3}) is inherently
sequential, thus it is difficult to parallelize it over time. This implies that we need to solve the sequence of linear systems
(\ref{eq1.9a}) one by one. Then Krylov subspace methods with suitable preconditioners \cite{Lei2013x,Pan2014,Donatelli16}
can be efficient candidates for solving Toeplitz-like linear systems since their complexity is of only  $\mathcal{O}(N\log N)$ arithmetic operations per iteration step.

In order to solve Eq. (\ref{eq1.9a}) effectively, we consider two specific classes of problems:
\begin{itemize}
\item[\romannumeral1)] When the diffusion coefficient $\xi(x,t) \equiv \xi$, the coefficient matrix of Eq. (\ref{eq1.9a})
will be a time-independent Toeplitz matrix, i.e. $\mathcal{M}^{j+1} = \mathcal{M}$; then we can compute its matrix inverse
via the Gohberg-Semencul formula (GSF) \cite{Gohberg92} using only its first and last columns. Such a strategy does not need
to call the preconditioned Krylov subspace solvers at each time level $0\leq j\leq M-1$, and the solution at each time
level (i.e., $\mathcal{M}^{-1}{\bm u}^{j+1}$) can be calculated via about six FFTs, thus saving considerable computational cost;
refer to \cite{Gu2017,Huang2017,Lin2018,Zhao2019} for detail.
\item[\romannumeral2)] When the diffusion coefficient is just a function related to both $x$ and $t$, i.e., $\xi(x,t)$, the coefficient matrix of Eq.
(\ref{eq1.9a}) becomes the sum of a scalar matrix and of a diagonal-multiply-Toeplitz matrix, which is time-dependent. In this case,
Eq. (\ref{eq1.9a}) has to be solved via a preconditioned Krylov subspace solver at each time level $j$.
\end{itemize}
Based on the above considerations, we still require to solve several nonsymmetric Toeplitz(-like) linear systems, whose matrix-vector
products can be efficiently calculated via FFTs, thus we utilize the biconjugate gradient stabilized (BiCGSTAB) method which has a fast and smooth convergence \cite{Vorst92}. For accelerating BiCGSTAB, we consider the following skew-circulant and banded
preconditioners:
\begin{equation}
P_{sk} = \begin{cases}
c_0I - \frac{\xi}{h^{\alpha}}\left[p\cdot sk(W_{\alpha}) + (1 - p)sk(W^{T}_{\alpha})\right],& \xi(x,t) \equiv \xi,\\
c_0I - \frac{\xi^{(j+1)}}{h^{\alpha}}\left[p\cdot sk(W_{\alpha}) + (1 - p)sk(W^{T}_{\alpha})\right],& \xi^{(j+1)} =
\frac{1}{N-1}\sum\limits^{N-1}_{i=1}\xi(x_i,t_{j+1}),
\end{cases}
\label{eq3.4}
\end{equation}
where the vector ${\bm \delta} = [w^{(\alpha)}_1,w^{(\alpha)}_2,\cdots,w^{(\alpha)}_{N-2},-w^{(\alpha)}_0]^T$ is the first column of
the skew-circulant matrix $sk(W_{\alpha})$ \cite{Zhao2019}, and
\begin{equation}
P_b =
\begin{cases}
c_0I - \frac{\xi}{h^{\alpha}}\left[pW_{\alpha,\ell} + (1 - p)W^{T}_{\alpha,\ell}\right],&\xi(x,t) \equiv \xi,\\
c_0I - \frac{K^{(j+1)}}{h^{\alpha}}\left[pW_{\alpha,\ell} + (1 - p)W^{T}_{\alpha,\ell}\right],&{\rm(general~case)},
\end{cases}
\label{eq3.5}
\end{equation}
with the band matrix
\begin{equation*}
W_{\alpha,\ell} =
\begin{bmatrix}
w^{(\alpha)}_1           & w^{(\alpha)}_0 \\
\vdots                   & w^{(\alpha)}_1 & w^{(\alpha)}_0   \\
w^{(\alpha)}_{\ell} &                & \ddots & \ddots\\
                         & \ddots         &                     &\ddots  &w^{(\alpha)}_0\\%
                         &                & w^{(\alpha)}_{\ell} & \cdots & w^{(\alpha)}_1 \\
\end{bmatrix},\quad \ell\in \mathbb{N}^{+},
\end{equation*}
respectively. Meanwhile, the high efficiency of skew-circulant and banded preconditioners for (time-)space FDEs has been shown in
\cite{Lin2014,Zhao16,Zhao2019}.

In practical implementations, when $P_{sk}$ or $P_b$ is employed as the preconditioner, a fast preconditioned
version of the BiCGSTAB method is obtained. During each BiCGSTAB iteration, two preconditioning steps are added in which one has to solve either the linear system $P_{sk}{\bm z} = {\bm y}$ or $P_b {\bm z} = {\bm y}$ for some given vector ${\bm y}$. Thus, some additional storage and computational cost are still required. However, we point out that $P_{sk}$ (resp., $P_b$) can also be efficiently
stored in $\mathcal{O}(N)$ (resp., $\mathcal{O}(\ell N)$) memory by only storing the $(N-1)$-dimensional vector ${\bm \delta}$
in (\ref{eq3.4}) (resp., the band matrix $W_{\alpha,\ell}$ in (\ref{eq3.5})). Besides, as $P_{sk}$ is the skew matrix\footnote{
If the diffusion coefficient $\xi(x,t) \equiv \xi(t)$, then $\xi^{(j+1)}$ are time-varying constants, which is available
for other similar cases.}, we observe
that
\begin{equation}
P_{sk} = \Omega^{*}F^{*}\Big\{c_0I - \frac{\xi^{(j+1)}}{h^{\alpha}}\left[p\Lambda_s + (1 - p)\bar{\Lambda}_s)\right]\Big\}F\Omega,\quad
sk(W_{\alpha}) = \Omega^{*}F^{*}\Lambda_s F\Omega,
\label{eq3.6}
\end{equation}
where $\Omega = {\rm diag}\left(1,(-1)^{-\frac{-1}{N-1}},\cdots,(-1)^{-\frac{N-2}{N-1}}\right)$, $F$ is the discrete
Fourier matrix and its conjugate transpose $F^{*}$. According to Eq. (\ref{eq3.6}), the inverse-matrix-vector product
${\bm z} = P^{-1}_{sk}{\bm y}$ can be carried out in $\mathcal{O}(N \log N)$ operations via the (inverse) FFTs. Most
importantly, the diagonal matrix $\Lambda_s$ can be computed in advance and only once per time step. On
the other hand, since $W_{\alpha,\ell}$ is a band matrix, then $P_b$ should be a band matrix of bandwidth $2\ell + 1$
and ${\bm z} = P^{-1}_b{\bm y}$ can be computed by the banded LU decomposition \cite{Lin2014,Zhao16} in $\mathcal{O}(\ell
N)$ arithmetic operations ($\ell \ll N$). In one word, we employ a fast preconditioned BiCGSTAB solution method with
low memory requirement and computational cost per iteration, while the number of iterations and thus the total computational
cost are greatly reduced. Compared to the skew-circulant preconditioner, the banded preconditioner needs more computational
cost to update at each time level; refer to the next section for a discussion.

On the other hand, it is worthwhile to note that when $\alpha\geq \alpha_0$, the coefficient matrix $\mathcal{M}^{(j+1)}$
are diagonally dominant with positive diagonal elements \cite{Lin2020} due to Proposition \ref{rem2.1} and $\xi(x,t) > 0$. Meanwhile, the
banded preconditioner was shown to be considerably efficient for solving the linear systems with diagonally dominant coefficient
matrix, which arise from the numerical discretization of (time-)space FDEs; refer, e.g., to \cite{Lin2014,Zhao16,Lin2020} for
a discussion.

\section{Numerical experiments}
\label{sec4}
The numerical experiments presented in this section have a two-fold objective. They illustrate that the
proposed implicit difference scheme (IDS) for the GTSFDE (\ref{eq1.1}) can indeed converge with
the order of $\mathcal{O}(\tau^{2 - \gamma} + h^2)$. Meanwhile, they assess the computational
efficiency of the fast solution techniques described in Section \ref{sec3}. Our choice of Krylov subspace method and direct solver in
Example 2
(where $\xi(x,t)\equiv \xi$ and $\mathcal{M}^{(j+1)}\equiv\mathcal{M}$ will be independent of time levels)
are the built-in MATLAB implementations of the preconditioned BiCGSTAB method and of the LU factorization, respectively, while in Example 1 with variable coefficients (where the coefficient matrices $\mathcal{M
}^{(j+1)}$ change at each time level) we use the MATLAB's backslash operator. The stopping criterion for the BiCGSTAB method with the two different preconditioners is $\|{\bm r}^{(k)}\|_{2}/\|{\bm r}^{(0)}\|_2 \leq 10^{-12}$, where ${\bm r}^{(k)}$ is the residual vector of the linear system after $k$ iterations; the initial guess is chosen as the zero vector. All experiments were performed on a Windows 10 (64 bit) PC-Intel(R) Core(TM) i5-8250U CPU @1.60 GHz--1.80GHz, 8 GB of RAM using MATLAB 2017b with machine epsilon $10^{-16}$ in double precision floating point arithmetic. By the way, all timings (measured in seconds) are averages calculated over 20 runs of our algorithms. Before we report the numerical results of the IDS for the problem (\ref{eq1.1}), we introduce the following notations that are adopted throughout this section:
\[
{\rm Error}_{\infty} = \max_{0\leq j\leq M}\|{\bm E}^{j}\|_{\infty}\quad {\rm and}\quad {\rm Error}_2 =
\max_{0\leq j\leq M}\|{\bm E}^{j}\|_2,
\]
then
\begin{equation*}
{\rm Rate}_{\infty} =
\begin{cases}
\log_{\tau_1/\tau_2}\left(\frac{{\rm Error}_{\infty,\tau_1}}{{\rm Error}_{\infty,\tau_2}}\right),&{\rm(temporal~convergence~order)},\\
\log_{h_1/h_2}\left(\frac{{\rm Error}_{\infty,h_1}}{{\rm Error}_{\infty,h_2}}\right),&{\rm(spatial~convergence~order)},
\end{cases}
\end{equation*}
and
\begin{equation*}
{\rm Rate}_2 =
\begin{cases}
\log_{\tau_1/\tau_2}\left(\frac{{\rm Error}_{2,\tau_1}}{{\rm Error}_{2,\tau_2}}\right),&{\rm(temporal~convergence~order)},\\
\log_{h_1/h_2}\left(\frac{{\rm Error}_{2,h_1}}{{\rm Error}_{2,h_2}}\right),&{\rm(spatial~convergence~order)}.
\end{cases}
\end{equation*}
\noindent\textbf{Example 1.} In this example, we solve the initial-boundary value problem of GTSFDE (\ref{eq1.1})
with variable coefficients and $\lambda(t) = e^{-bt}, b \geq 0$, the spatial domain $[x_L,x_R] = [0,2]$ and the
time interval is $[0,T] = [0,1]$. The diffusion coefficient function is given as $\xi(x,t) = 1 + x^2 + \sin t$.
The source term is
\begin{equation*}
\begin{split}
f(x,t) = & \frac{2t^{3 - \gamma}e^{-bt}}{\Gamma(4 - \gamma)}x^2(2-x)^2 - g(t)\xi(x,t)
\Bigg\{\frac{4\Gamma(3)}{\Gamma(3 - \alpha)}\Big[px^{2-\alpha} + (1 - p)(2-x)^{2-\alpha}
\Big]\\
& - \frac{4\Gamma(4)}{\Gamma(4 - \alpha)}\Big[px^{3-\alpha} + (1 - p)(2-x)^{3-\alpha}
\Big] + \frac{\Gamma(5)}{\Gamma(5-\alpha)}\Big[px^{4-\alpha} + (1 - p)(2-x)^{4-\alpha}\Big]\Bigg\},
\end{split}
\end{equation*}
and the initial-boundary value conditions are
\[u(x,0) = g(0)x^2(2-x)^2,\qquad \mathrm{and}\qquad u(0,t) = u(2,t)=0.\]
The exact (and smooth) solution to this problem is $u(x,t) = g(t)x^2(2-x)^2$,
where $g(t)$ is given as follows:
\begin{equation*}
g(t) = 1 + \frac{2 - (2 + 2bt + b^2t^2)e^{-bt}}{b^3}
\end{equation*}
for any $\alpha\in(1,2)$ and $b\in \mathbb{R}^{+}$. Numerical experiments with our proposed difference scheme are reported in the following Tables \ref{tab5}--\ref{tab8}.

\begin{table}[!htpb]
\caption{$L_2$-norm and maximum norm errors versus grid size reduction when
$h = 2^{-12}$ and $p = 0.7$ in Example 1.}
\centering
\begin{tabular}{crcccccccc}
\hline & \multicolumn{4}{c}{$b = 1.0$} &\multicolumn{4}{c}{$b = 2.0$}\\
[-2pt]\cmidrule(r{0.7em}r{0.7em}){3-6}\cmidrule(r{0.7em}r{0.6em}){7-10}\\[-11pt]
$(\gamma,\alpha)$ &$\tau$ &Error$_{\infty}$ &Rate$_{\infty}$ &
Error$_2$ & Rate$_2$ &Error$_{\infty}$ &Rate$_{\infty}$ &Error$_2$ & Rate$_2$     \\
\hline
(0.2,1.1) &1/8  &5.9654e-4 &--     &5.5779e-4 &--     &3.1311e-4 &--     &2.9126e-4 &--     \\
          &1/16 &1.7385e-4 &1.7788 &1.6250e-4 &1.7793 &9.0388e-5 &1.7925 &8.4009e-5 &1.7937 \\
          &1/32 &5.0703e-5 &1.7777 &4.7379e-5 &1.7781 &2.6194e-5 &1.7869 &2.4335e-5 &1.7875 \\
          &1/64 &1.4813e-5 &1.7752 &1.3843e-5 &1.7751 &7.6240e-6 &1.7806 &7.0838e-6 &1.7804 \\
\hline
(0.5,1.5) &1/8  &1.0328e-3 &--     &1.0162e-3 &--     &5.1328e-4 &--     &5.0407e-4 &--     \\
          &1/16 &3.7458e-4 &1.4632 &3.6869e-4 &1.4627 &1.8639e-4 &1.4614 &1.8284e-4 &1.4630 \\
          &1/32 &1.3450e-4 &1.4777 &1.3235e-4 &1.4781 &6.7060e-5 &1.4748 &6.5809e-5 &1.4742 \\
          &1/64 &4.8098e-5 &1.4836 &4.7330e-5 &1.4835 &2.4016e-5 &1.4815 &2.3557e-5 &1.4821 \\
\hline
(0.9,1.9) &1/8  &2.9303e-3 &--     &2.8851e-3 &--     &1.3940e-3 &--     &1.3710e-3 &--     \\
          &1/16 &1.3909e-3 &1.0750 &1.3700e-3 &1.0744 &6.6816e-4 &1.0610 &6.5664e-4 &1.0621 \\
          &1/32 &6.5575e-4 &1.0848 &6.4585e-4 &1.0849 &3.1678e-4 &1.0767 &3.1149e-4 &1.0759 \\
          &1/64 &3.0744e-4 &1.0928 &3.0279e-4 &1.0929 &1.4894e-4 &1.0887 &1.4650e-4 &1.0883 \\
\bottomrule
\end{tabular}
\label{tab5}
\end{table}
\begin{table}[!htpb]
\caption{$L_2$-norm and maximum norm errors versus grid size reduction when
$\tau = 2^{-10}$ and $p = 0.7$ in Example 1.}
\centering
\begin{tabular}{crcccccccc}
\hline &\multicolumn{4}{c}{$b = 1.0$} &\multicolumn{4}{c}{$b = 2.0$}\\
[-2pt]\cmidrule(r{0.7em}r{0.7em}){3-6}\cmidrule(r{0.7em}r{0.6em}){7-10}\\[-11pt]
$(\gamma,\alpha)$ &$h$ &Error$_{\infty}$ &Rate$_{\infty}$ &
Error$_2$ & Rate$_2$ &Error$_{\infty}$ &Rate$_{\infty}$ &Error$_2$ & Rate$_2$     \\
\hline
(0.2,1.1) &2/8  &1.0332e-1 &--     &9.5781e-2   &--     &1.0756e-1 &--     &9.6687e-2 &--     \\
          &2/16 &2.4194e-2 &2.0944 &2.3302e-2   &2.0393 &2.3916e-2 &2.1691 &2.3510e-2 &2.0400 \\
          &2/32 &7.3546e-3 &1.7180 &5.5686e-3   &2.0650 &6.6797e-3 &1.8401 &5.6175e-3 &2.0653 \\
          &2/64 &2.0330e-3 &1.8550 &1.3355e-3   &2.0599 &1.8355e-3 &1.8636 &1.3477e-3 &2.0594 \\
\hline
(0.5,1.5) &2/8  &7.0414e-2 &--     &6.7030e-2   &--     &6.9027e-2 &--     &6.5647e-2 &--     \\
          &2/16 &1.6525e-2 &2.0912 &1.5689e-2   &2.0951 &1.6114e-2 &2.0988 &1.5317e-2 &2.0996 \\
          &2/32 &3.9248e-3 &2.0740 &3.7129e-3   &2.0791 &3.8292e-3 &2.0732 &3.6158e-3 &2.0827 \\
          &2/64 &1.0322e-3 &1.9269 &8.8843e-4   &2.0632 &9.5842e-4 &1.9983 &8.6283e-4 &2.0672 \\
\hline
(0.9,1.9) &2/8  &6.9963e-2 &--     &7.0620e-2   &--     &6.6930e-2 &--     &6.7553e-2 &--     \\
          &2/16 &1.7061e-2 &2.0359 &1.7145e-2   &2.0423 &1.6307e-2 &2.0372 &1.6387e-2 &2.0435 \\
          &2/32 &4.1828e-3 &2.0282 &4.1803e-3   &2.0361 &3.9886e-3 &2.0315 &3.9871e-3 &2.0391 \\
          &2/64 &1.0354e-3 &2.0143 &1.0281e-3   &2.0236 &9.7927e-4 &2.0261 &9.7271e-4 &2.0353 \\
\bottomrule
\end{tabular}
\label{tab6}
\end{table}
In Tables \ref{tab5}--\ref{tab6}, we display the maximum-norm errors and $L_2$-norm errors of the
IDS scheme for solving the problem (\ref{eq1.1}) with variable diffusion coefficients in spatial and temporal
variables, respectively. More precisely, the results of Table \ref{tab5} with different $(\gamma,\alpha,b)$'s
and $h = 2^{-12}$ show, as expected, a reduction in the maximum- or $L_2$-norm error as the number of time steps of our IDS is increased, and the temporal convergence order of IDS is $\mathcal{O}(\tau^{2 -
\gamma})$. At the same time, our experiments displayed in Table \ref{tab6} with different $(\gamma,\alpha,b)$'s and $\tau = 2^{-10}$ show a reduction in the maximum- or $L_2$-norm error as the size of time steps of our IDS is decreased, and thus the convergence order in space is $\mathcal{O}(h^2)$. In conclusion, the
numerical convergence orders are consistent with the theoretical estimate $\mathcal{O}(\tau^{2 - \gamma} +
h^2)$ presented in Section \ref{sec2.2}.

\begin{table}[!htpb]
\caption{Numerical comparisons of the direct, iterative, and preconditioned iterative methods
for solving Example 1 with $\tau = 2^{-12}$, $b=1.0$, and $p = 0.7$.}
\centering
\begin{tabular}{lrrcccccc}
\hline &&Direct&\multicolumn{2}{c}{Noprec} &\multicolumn{2}{c}{Banded($\ell = 8$)}&\multicolumn{2}{c}{Skew-cir}\\
[-2pt]\cmidrule(r{0.7em}r{0.7em}){4-5}\cmidrule(r{0.7em}r{0.6em}){6-7}\cmidrule(r{0.7em}r{0.6em}){8-9}\\[-11pt]
$(\gamma,\alpha)$ &$N$ &CPU(s) &Iter &CPU(s) &Iter &CPU(s) &Iter &CPU(s)  \\
\hline
(0.2,1.1) 
          &128  &25.895   &225.0 &84.890  &6.0  &23.462  &13.8 &26.027   \\
          &256  &55.777   &994.5 &893.026 &7.3  &32.152  &14.3 &31.613   \\
          &512  &249.966  &\dag  &\dag    &9.9  &54.337  &15.0 &56.042   \\
          &1024 &1859.940 &\dag  &\dag    &16.2 &102.047 &15.6 &75.748   \\
\hline
(0.5,1.5) 
          &128  &25.666   &65.0  &38.702   &5.5  &23.111  &14.1 &26.180  \\
          &256  &52.250   &127.6 &68.884   &8.0  &32.794  &14.8 &31.585  \\
          &512  &249.797  &261.8 &311.449  &11.7 &57.665  &15.6 &57.493  \\
          &1024 &1836.925 &999.9 &3358.602 &18.7 &109.722 &16.3 &77.504  \\
\hline
(0.9,1.9) 
          &128  &24.118   &39.0  &29.815   &3.0  &21.174  &12.0 &24.843  \\
          &256  &52.021   &71.8  &49.841   &4.0  &26.891  &13.4 &30.941  \\
          &512  &260.412  &147.4 &183.047  &5.0  &42.013  &14.8 &55.654  \\
          &1024 &1835.804 &301.3 &421.588  &7.0  &70.558  &15.7 &76.047  \\
\bottomrule
\end{tabular}
\label{tab7}
\end{table}
\begin{table}[!htpb]
\caption{Numerical comparisons of the direct, iterative, and preconditioned iterative methods
for solving Example 1 with $\tau = 2^{-12}$, $b = 2.0$ and $p = 0.3$.}
\centering
\begin{tabular}{lrrcccccc}
\hline &&Direct&\multicolumn{2}{c}{Noprec} &\multicolumn{2}{c}{Banded($\ell = 8$)}&\multicolumn{2}{c}{Skew-cir}\\
[-2pt]\cmidrule(r{0.7em}r{0.7em}){4-5}\cmidrule(r{0.7em}r{0.6em}){6-7}
\cmidrule(r{0.7em}r{0.6em}){8-9}\\[-11pt]
$(\gamma,\alpha)$ &$N$ &CPU &Iter &CPU &Iter &CPU &Iter &CPU  \\
\hline
(0.2,1.1) 
          &128  &26.045   &204.3 &80.795   &5.5  &23.267  &15.9 &26.447  \\
          &256  &54.872   &981.7 &1017.883 &7.3  &30.687  &15.9 &32.356  \\
          &512  &250.823  &\dag  &\dag     &11.0 &57.626  &15.8 &57.388  \\
          &1024 &1822.899 &\dag  &\dag     &17.8 &107.799 &16.6 &77.719  \\
\hline
(0.5,1.5) 
          &128  &26.470   &60.0  &36.100   &5.4  &23.145  &14.2 &26.188  \\
          &256  &54.736   &123.1 &67.480   &7.8  &33.156  &14.8 &31.489  \\
          &512  &260.543  &246.8 &293.753  &11.3 &57.874  &14.9 &55.981  \\
          &1024 &1824.693 &502.1 &658.004  &18.9 &110.868 &15.8 &76.082  \\
\hline
(0.9,1.9) 
          &128  &24.914   &39.4  &30.340   &3.0  &22.314  &11.9 &24.145  \\
          &256  &52.028   &70.3  &47.876   &3.0  &27.789  &13.6 &30.749  \\
          &512  &250.682  &144.2 &175.927  &5.0  &42.018  &14.9 &55.671  \\
          &1024 &1823.420 &298.9 &410.570  &7.0  &70.534  &15.9 &75.984  \\
\bottomrule
\end{tabular}
\label{tab8}
\end{table}
In Tables \ref{tab7}--\ref{tab8}, the performance of the direct, iterative, and preconditioned
iterative methods for Eq.~(\ref{eq1.9a}) are illustrated along with the elapsed CPU time and the average number of iterations. Here the symbols ``Direct", ``Noprec", ``Banded($\ell = 8$)" and ``Skew-cir" mean that the
sequence of linear systems (\ref{eq1.9a}) is consecutively solved by using the MATLAB's backslash operator or the BiCGSTAB routine available in MATLAB with no preconditioner, banded preconditioner and skew-circulant preconditioner, respectively. As can be seen from Tables \ref{tab7}--\ref{tab8}, both banded and skew-circulant preconditioners are fairly efficient to accelerate the BiCGSTAB method for solving Eq. (\ref{eq1.9a}) in terms of the elapsed CPU time and the number of iterations, especially when the number of grid nodes increases. Moreover, we remark that $P_{sk}$ exhibits more robust performance than $P_b$ in terms of average number of iterations, i.e., compared to BiCGSTAB with $P_b$, and the average of number of BiCGSTAB with $P_{sk}$ is weakly sensitive to the spatial grid size. In addition, it can be observed that the performance of BiCGSTAB with $P_b$ becomes better when $\alpha = 1.9$, because the banded preconditioner has been proved to be very efficient for solving Eq.~(\ref{eq1.9a}), whose coefficient matrices are diagonally dominant -- cf. Proposition \ref{rem2.1} and Section \ref{sec3} for a discussion. In conclusion, the skew-circulant preconditioner is still recommended for enhancing the convergence of BiCGSTAB applied to solve Eq. (\ref{eq1.9a}), when $1 < \alpha < \alpha_0$, whereas the banded preconditioner is recommended if $\alpha \geq \alpha_0$, because the coefficient matrices are diagonally dominant -- cf. Section \ref{sec3}.
\vspace{1mm}

\noindent\textbf{Example 2}. In this numerical example, we show the effect on the solution due to the presence of a singularity in the temporal derivative. For clarity, we assume that we can isolate a nonsmooth part from $u(x,t)$ as is shown in \cite{Alikhanov17}. We suppose that the solution to the problem~(\ref{eq1.1})
with $f =0,~u(0,t)=u(1,t)=0,~u(x,0) = 5x^3(1 - x)^3$, $\lambda(t) = e^{-bt}$, $\xi(x,t) \equiv
\xi$ and $[x_L,x_R]\times[0,T] = [0,1]^2$ has the following form
\begin{equation}
u(x,t) = 5x^3(1 - x)^3\left[1 - \frac{\sqrt{t}}{\Gamma(1.5)}e^{-bt}\right] + v(x,t),
\label{eq4.1}
\end{equation}
where $v(x,t)$ is the exact solution to the problem
\begin{equation}
\begin{split}
{}^{C}_{0}D^{\gamma,\lambda(t)}_{t}v(x,t) & = \xi\Big[p{}_{a}D^{\alpha}_{x}v(x,t) +
(1 - p){}_xD^{\alpha}_bv(x,t)\Big] + 5x^3(1 - x)^3
\left[\frac{t^{\frac{1}{2} - \gamma}}{\Gamma(\frac{3}{2} - \gamma)} - \frac{bt^{\frac{3}{2}- \gamma}}{\Gamma(\frac{5}{2} - \gamma)}\right]e^{-bt}
\\
&\quad~ + 5\xi \cdot q(x)\left[1 - \frac{\sqrt{t}}{\Gamma(1.5)}e^{-bt}\right],\quad (x,t)\in(0,1)\times(0,T],
\end{split}
\end{equation}
along with the initial condition $v(x,0) = 0$~($x\in[0,1]$), the boundary conditions $
v(0,t) = v(1,t) = 0$~($t\in(0,T]$), and with
\begin{equation}
\begin{split}
q(x) & = \frac{\Gamma(4)}{\Gamma(4- \alpha)}\Big[px^{3-\alpha} + (1-p)(1-x)^{3 -
\alpha}\Big] - \frac{3\Gamma(5)}{\Gamma(5- \alpha)}\Big[px^{4-\alpha} + (1-p)(1-x)^{4 -
\alpha}\Big]~ + \\
&\quad~\frac{3\Gamma(6)}{\Gamma(6 - \alpha)}\Big[px^{5-\alpha} + (1-p)(1-x)^{5-\alpha}
\Big] - \frac{\Gamma(7)}{\Gamma(7 -\alpha)}\Big[px^{6-\alpha} + (1-p)(1-x)^{6-\alpha}\Big].
\end{split}
\end{equation}
Since $u(x,t)$ is not available in this case, we estimate both the temporal and the spatial convergence orders of our difference scheme by computing the approximate solution on two different grids $\varpi_{h,\tau_1}$ and $\varpi_{h,\tau_2}$, where $\varpi_{h,\tau_1}\subset \varpi_{h,\tau_2}$. The numerical calculations of $u(x,t)$ by formula (\ref{eq4.1})
are reported in Tables \ref{tab9}--\ref{tab12}, where we focus on the convergence order of numerical solutions at
the \textit{final time point}\footnote{It is worth noting that the results of our experiments
showing that the temporal convergence order of our proposed scheme in the \textit{only} $L^2$-norm (which is the same as Table \ref{tab9}) in the whole domain is 1 but less than $2-\gamma$ still remains consistent with the stability analysis and error estimates based on the $L^2$-norm introduced in Section \ref{sec2.2}.}.
\begin{table}[!htpb]
\caption{$L_2$-norm and maximum norm errors versus grid size reduction when
$h = 2^{-11}$ and $p = 0.4$ in Example 2.}
\centering
\begin{tabular}{crrcccccccc}
\hline & & & \multicolumn{4}{c}{$b = 3.0$} &\multicolumn{4}{c}{$b = 4.0$}\\
[-2pt]\cmidrule(r{0.7em}r{0.7em}){4-7}\cmidrule(r{0.7em}r{0.6em}){8-11}\\[-11pt]
$(\gamma,\alpha)$ &$\tau_1$ &$\tau_2$ &Error$_{\infty}$ &Rate$_{\infty}$ &
Error$_2$ & Rate$_2$ &Error$_{\infty}$ &Rate$_{\infty}$ &Error$_2$ & Rate$_2$     \\
\hline
(0.2,1.2) &1/20  &1/40  &2.344e-5 &--    &1.637e-5 &--    &1.451e-5 &--    &1.019e-5 &--    \\
          &1/40  &1/80  &1.185e-5 &0.984 &8.280e-6 &0.983 &7.142e-6 &1.022 &5.021e-6 &1.021 \\
          &1/80  &1/160 &6.016e-6 &0.978 &4.206e-6 &0.977 &3.574e-6 &0.999 &2.514e-6 &0.998 \\
          &1/160 &1/320 &3.054e-6 &0.978 &2.136e-6 &0.978 &1.800e-6 &0.990 &1.266e-6 &0.989 \\
\hline
(0.5,1.5) &1/20  &1/40  &7.456e-6 &--    &5.225e-6 &--    &4.672e-6 &--    &3.283e-6 &--    \\
          &1/40  &1/80  &3.452e-6 &1.111 &2.419e-6 &1.111 &2.066e-6 &1.178 &1.452e-6 &1.177 \\
          &1/80  &1/160 &1.633e-6 &1.079 &1.145e-6 &1.079 &9.438e-7 &1.130 &6.634e-7 &1.130 \\
          &1/160 &1/320 &7.850e-7 &1.057 &5.502e-7 &1.057 &4.418e-7 &1.095 &3.106e-7 &1.095 \\
\hline
(0.9,1.9) &1/20  &1/40  &4.812e-6 &--    &3.322e-6 &--    &3.696e-6 &--    &2.554e-6 &--    \\
          &1/40  &1/80  &2.279e-6 &1.078 &1.574e-6 &1.078 &1.727e-6 &1.098 &1.193e-6 &1.098 \\
          &1/80  &1/160 &1.074e-6 &1.086 &7.410e-7 &1.086 &8.074e-7 &1.097 &5.578e-7 &1.097 \\
          &1/160 &1/320 &5.036e-7 &1.092 &3.476e-7 &1.092 &3.773e-7 &1.098 &2.606e-7 &1.098 \\
\bottomrule
\end{tabular}
\label{tab9}
\end{table}
\begin{table}[!htpb]
\caption{$L_2$-norm and maximum norm errors versus grid size reduction when $\tau =
2^{-11}$ and $p = 0.4$ in Example 2.}
\centering
\begin{tabular}{crrcccccccc}
\hline & & & \multicolumn{4}{c}{$b = 3.0$} &\multicolumn{4}{c}{$b = 4.0$}\\
[-2pt]\cmidrule(r{0.7em}r{0.7em}){4-7}\cmidrule(r{0.7em}r{0.6em}){8-11}\\[-11pt]
$(\gamma,\alpha)$ &$h_1$ &$h_2$ &Error$_{\infty}$ &Rate$_{\infty}$ &
Error$_2$ & Rate$_2$ &Error$_{\infty}$ &Rate$_{\infty}$ &Error$_2$ & Rate$_2$     \\
\hline
(0.2,1.2) &1/10  &1/20  &2.318e-3 &--    &1.500e-3  &--    &2.424e-3 &--    &1.555e-3 &--    \\
          &1/20  &1/40  &5.285e-4 &2.133 &3.437e-4  &2.126 &5.517e-4 &2.135 &3.557e-4 &2.128 \\
          &1/40  &1/80  &1.295e-4 &2.029 &8.427e-5  &2.028 &1.350e-4 &2.031 &8.687e-5 &2.034 \\
          &1/80  &1/160 &5.271e-5 &1.297 &2.141e-5  &1.977 &3.369e-5 &2.003 &2.167e-5 &2.003 \\
\hline
(0.5,1.5) &1/10  &1/20  &2.003e-3 &--    &1.274e-3  &--    &2.084e-3 &--    &1.324e-3 &--    \\
          &1/20  &1/40  &4.577e-4 &2.130 &2.949e-4  &2.111 &4.760e-4 &2.130 &3.063e-4 &2.111 \\
          &1/40  &1/80  &1.122e-4 &2.028 &7.208e-5  &2.033 &1.167e-4 &2.029 &7.487e-5 &2.033 \\
          &1/80  &1/160 &2.799e-5 &2.003 &1.790e-5  &2.009 &2.908e-5 &2.004 &1.860e-5 &2.009 \\
\hline
(0.9,1.9) &1/10  &1/20  &1.445e-3 &--    &8.723e-4  &--    &1.501e-3 &--    &9.057e-4 &--    \\
          &1/20  &1/40  &3.256e-4 &2.150 &2.059e-4  &2.083 &3.382e-4 &2.150 &2.138e-4 &2.083 \\
          &1/40  &1/80  &7.940e-5 &2.036 &5.069e-5  &2.022 &8.245e-5 &2.036 &5.262e-5 &2.022 \\
          &1/80  &1/160 &1.975e-5 &2.007 &1.262e-5  &2.006 &2.051e-5 &2.007 &1.310e-5 &2.006 \\
\bottomrule
\end{tabular}
\label{tab10}
\end{table}
\begin{table}[!htpb]
\centering
\caption{Numerical comparisons of the direct, iterative, and preconditioned iterative methods
for solving Example 2 with $b = 3.0$, $\tau = 2^{-12}$, $\kappa = 5$, and $p = 0.4$.}
\begin{tabular}{lrrrrccccc}
\hline &&\multicolumn{2}{c}{Direct} &\multicolumn{3}{c}{Banded($\ell = 8$)}&\multicolumn{3}{c}{Skew-cir}\\
[-2pt]\cmidrule(r{0.7em}r{0.7em}){3-4}\cmidrule(r{0.7em}r{0.6em}){5-7}
\cmidrule(r{0.7em}r{0.6em}){8-10}\\[-11pt]
$(\gamma,\alpha)$ &$N$  &CPU    &Memory &Iter &CPU    &Memory &Iter &CPU    &Memory  \\
\hline
(0.2,1.2)         &128  &17.357 &13.003 &16.0 &17.458 &12.254 &15.0 &16.972 &12.148 \\
                  &256  &21.690 &27.631 &25.0 &21.218 &24.587 &15.0 &20.873 &24.170 \\
                  &512  &30.846 &62.137 &36.0 &29.534 &48.653 &16.0 &29.218 &48.215 \\
                  &1024 &50.219 &152.15 &54.0 &44.977 &97.184 &17.0 &44.552 &96.305 \\
\hline
(0.5,1.5)         &128  &17.339 &13.003 &15.0 &17.292 &12.254 &14.0 &16.941 &12.148 \\
                  &256  &21.646 &27.631 &24.0 &21.144 &24.587 &14.0 &20.852 &24.170 \\
                  &512  &31.001 &62.137 &38.0 &30.003 &48.653 &16.0 &29.225 &48.215 \\
                  &1024 &50.261 &152.15 &61.0 &45.465 &97.184 &16.0 &44.498 &96.305 \\
\hline
(0.9,1.9)         &128  &17.297 &13.003 &8.0  &16.919 &12.254 &11.0 &16.928 &12.148 \\
                  &256  &21.588 &27.631 &10.0 &20.928 &24.587 &12.0 &20.837 &24.170 \\
                  &512  &30.901 &62.137 &16.0 &29.307 &48.653 &12.0 &29.086 &48.215 \\
                  &1024 &50.304 &152.15 &27.0 &44.562 &97.184 &12.0 &44.401 &96.305 \\
\bottomrule
\end{tabular}
\label{tab11}
\end{table}
\begin{table}[!htpb]
\centering
\caption{Numerical comparisons of the direct, iterative, and preconditioned iterative methods
for solving Example 2 with $b = 4.0$, $\tau = 2^{-12}$, $\kappa = 5$, and $p = 0.6$.}
\begin{tabular}{lrrcrccccc}
\hline &&\multicolumn{2}{c}{Direct} &\multicolumn{3}{c}{Banded($\ell = 8$)}&\multicolumn{3}{c}{Skew-cir}\\
[-2pt]\cmidrule(r{0.7em}r{0.7em}){3-4}\cmidrule(r{0.7em}r{0.6em}){5-7}
\cmidrule(r{0.7em}r{0.6em}){8-10}\\[-11pt]
$(\gamma,\alpha)$ &$N$  &CPU    &Memory &Iter &CPU    &Memory &Iter &CPU    &Memory  \\
\hline
(0.2,1.2)         &128  &17.289 &13.003 &16.0 &17.449 &12.254 &15.0 &16.959 &12.148 \\
                  &256  &21.568 &27.631 &25.0 &21.207 &24.587 &16.0 &20.908 &24.170 \\
                  &512  &30.891 &62.137 &36.0 &29.492 &48.653 &16.0 &29.211 &48.215 \\
                  &1024 &50.346 &152.15 &54.0 &44.959 &97.184 &17.0 &44.546	&96.305 \\
\hline
(0.5,1.5)         &128  &17.304 &13.003 &15.0 &17.237 &12.254 &14.0 &16.934 &12.148 \\
                  &256  &21.596 &27.631 &24.0 &21.106 &24.587 &14.0 &20.842 &24.170 \\
                  &512  &30.789 &62.137 &38.0 &29.892 &48.653 &16.0 &29.197	&48.215 \\
                  &1024 &50.273 &152.15 &61.0 &45.501 &97.184 &16.0 &44.501 &96.305 \\
\hline
(0.9,1.9)         &128  &17.297 &13.003 &8.0  &16.923 &12.254 &10.0 &16.915 &12.148 \\
                  &256  &21.583 &27.631 &10.0 &20.967 &24.587 &12.0 &20.852	&24.170 \\
                  &512  &30.841 &62.137 &16.0 &29.315 &97.184 &12.0 &29.101 &48.215 \\
                  &1024 &50.335 &152.15 &27.0 &44.496 &97.184 &12.0 &44.398 &96.305 \\
%
\bottomrule
\end{tabular}
\label{tab12}
\end{table}

As is seen from Table \ref{tab9}, in this example the temporal convergence order of our proposed difference scheme is almost 1 but is smaller than the theoretically estimated order --
$(2-\gamma)$, except for the values $(\gamma,\alpha) = (0.9,1.9)$. This can be explained by observing that
we selected only the part of the solution that yields a singularity in the first derivative of $u(x, t)$. This occurs when $v_t(x, t)$ is continuous and $v_{tt}(x, t)$ has a singularity at the initial point $
t = 0$. Moreover, Table \ref{tab10} shows that the spatial convergence order can still reach the theoretical estimate $\mathcal{O}(h^2)$, especially in the $L^2$ norm, when both the solution and the initial data are sufficiently smooth in the space variable $x$. In addition, the convergence order will be slightly better if $b$ increases, because the larger value of $b$ will make the solution $u(x,t)$ behave more smoothly. We conclude that in cases when we do not have enough information on the smoothness of the solution, we can calculate the convergence order as suggested above. If the estimated convergence order is smaller than
$\mathcal{O}(\tau^{2-\gamma} + h^2)$, then we should represent the solution as the sum of two functions,
one of which is non-smooth while the other is smooth but unknown~\cite{Alikhanov17}. Of course, finding such a representation may be difficult for the problem (\ref{eq1.1}) with variable coefficients, and this is an interesting research issue in its own right.

Tables \ref{tab11}--\ref{tab12} show the elapsed CPU time, the number of iterations and the memory cost (measured in megabytes) required by different solution techniques for solving the class \romannumeral1) of problem (\ref{eq1.1}) described in Section \ref{sec3}. Overall, the performance of the BiCGSTAB method with skew-circulant preconditioner are the best in terms of elapsed CPU time and memory cost, the banded preconditioner being a good alternative. However, the direct method (based on only one LU decomposition) is noncompetitive due to the large elapsed CPU time and memory costs, especially for the fine discretized meshes, because the large dense matrix and its LU decomposition factors need to be explicitly stored.
%
%
%
%
%
%
%

In addition, according to numerical results of Examples 1-2, it is interesting to observe that, although
we employ the fast preconditioned BiCGSTAB method to solve Eq. (\ref{eq1.9a}) corresponding to the IDS
(\ref{eq2.3}), the total CPU time is still high. In fact, the solution time comes from two main computations: 1) solving the sequence of linear systems (\ref{eq1.9a}); 2) evaluating the right-hand side vector of (\ref{eq1.9a}) by repeatedly summing the solutions at previous time levels. Our preconditioned
BiCGSTAB method can only alleviate the first cost, while we should analyse further the degradation of
CPU time due to handling the nonlocal property of the discrete temporal fractional derivative. However, such analysis is always difficult in the general case of $\lambda(t)$. In particular, if we set $\lambda(t) = e^{-bt}$ like in Examples 1-2, we can further alleviate the computational and memory cost of the proposed IDS. The derivation of such a more cost-effective scheme (\ref{eqA.5}) is presented in \ref{sec.esss}.
\vspace{1mm}
%
\section{Conclusions}
\label{sec5}
In this paper, the stability and convergence of an IDS scheme for solving the GTSFDEs with variable coefficients are studied via the diagonally weighted energy norm analysis. The proposed IDS can be proved to reach second order convergence in space and $(2 - \gamma)$-th approximation order in time for the GTSFDEs with variable coefficients. Moreover, numerical experiments involving problem (\ref{eq1.1}) with non-smooth solution are carried out yielding results completely in line with our theoretical analysis. The method can be easily extended to solve the variable coefficient GTSFDEs with other boundary conditions. Although the focus of the paper is on the case of one-dimensional spatial domains, the results can be extended to two- and three-dimensional domain; refer, e.g., to \cite{Vong2016}.

In addition, we have also shown an efficient implementation of the proposed IDS based on preconditioned iterative
solvers, achieving about $\mathcal{O}(N\log N)$ computational complexity and $\mathcal{O}(N)$ storage cost. Numerical evidence of the efficiency of the proposed preconditioning methods is reported. For the special choice of $\lambda(t) = e^{-bt}$, the fast sum-of-exponential approximations of the kernel in (\ref{eq1.2}) can be used to derive a cost-effective version of IDS (\ref{eqA.5}); then, numerical experiments are illustrated to show that the rate of the truncation error of this new IDS is about $\mathcal{O}(\tau^{2 - \gamma} + h^2)$. However, its rigorous stability and convergence analyses remain an open question. Meanwhile, numerical results show the fast IDS (\ref{eqA.5}) requires less CPU time and memory cost than the proposed IDS (\ref{eq2.3})

\appendix
\setcounter{lemma}{0}
\renewcommand{\thelemma}{\Alph{section}\arabic{lemma}}
\section{Fast SOE approximation of the generalized Caputo fractional derivative}
\label{sec.esss}
Due to the nonlocality of the generalized Caputo fractional derivative (\ref{eq1.2}), the proposed
scheme (\ref{eq1.9a}) requires the storage of the solution at all previous time steps which leads
to huge computational cost. This phenomenon also can be observed from the numerical experiments reported in
Section \ref{sec4}. To reduce the computational cost, we follow the work about fast \textit{L}1
formula \cite{Jiang17} for developing the SOE approximation of the generalized
Caputo fractional derivative with $\lambda(t) = e^{-bt}$, which is adopted in Section \ref{sec4}.
More precisely,
\begin{equation*}
\begin{split}
{}^{C}_0D^{\gamma,\lambda(t)}_tu(t)\mid_{t = t_j} & =
\frac{1}{\Gamma(1 - \gamma)}\int^{t_j}_0
\frac{e^{-b(t_j - s)}u'(s)ds}{(t_j - s)^{\gamma}}\\
& =
\frac{1}{\Gamma(1 - \gamma)}\int^{t_j}_{t_{j-1}}\frac{e^{-b(t_j - s)}u'(s)ds}{(t_j - s)^{\gamma}}
+ \frac{1}{\Gamma(1 - \gamma)}\int^{t_{j-1}}_0\frac{e^{-b(t_j - s)}u'(s)ds}{(t_j - s)^{\gamma}}\\
& = C_l(t_j) + C_h(t_j),
\end{split}
\end{equation*}
where the last equality defines the local part and the history part, respectively. For the local part,
we employ the generalized \textit{L}1 approximation recalled in Section \ref{sec2.1}, which approximates $u(s)$
on $[t_{j - 1},t_j]$ via a linear polynomial (with $u(t_{j-1})$ and $u(t_j)$ as the interpolation nodes)
or $u'(s)$ via a constant $\frac{u(t_j) - u(t_{j - 1})}{\tau}$. We have
\begin{equation}
\begin{split}
C_l(t_j) &\approx \frac{u(t_j) - u(t_{j - 1})}{\tau\Gamma(1 - \gamma)}\int^{t_j}_{t_{j-1}}
\frac{e^{-b(t_j - s)}ds}{(t_j - s)^{\gamma}}\\
& = \frac{u(t_j) - u(t_{j - 1})}{\tau\Gamma(2 - \gamma)}\Big(e^{-b\tau}\tau^{1 - \gamma} +
b\int^{\tau}_0e^{-b\theta}\theta^{1-\gamma}d\theta\Big),
\end{split}
\label{eqA.1}
\end{equation}
where the second integral can be evaluated via the MATLAB built-in function `\texttt{integral.m}'.
For the history part, we first recall the following lemma \cite{Jiang17} to approximate the history
part $C_h(t_j)$.
\begin{lemma}
Let $\epsilon$ denote tolerance error, $\delta$ cut-off time restriction and $T$ final time. Then
there are a natural number $N_{exp}$ and positive numbers $s_k$ and $w_k,~k = 1,2,\cdots,N_{exp}$
such that
\begin{equation*}
\left|\frac{1}{t^{\gamma}} - \sum^{N_{exp}}_{k=1}\omega_k e^{-s_kt}\right| < \epsilon,\quad
{\rm for~any~}t\in[\delta,T],
\end{equation*}
where $N_{exp} = \mathcal{O}((\log \epsilon^{-1})(\log\log \epsilon^{-1} + \log(T\delta^{-1}))
+ (\log\delta^{-1})(\log\log\epsilon^{-1} + \log\delta^{-1}))$.
\label{lemmaxa}
\end{lemma}

Therefore, when we set $\delta = \tau$ and apply Lemma \ref{lemmaxa}, then we obtain
\begin{equation}
\begin{split}
C_h(t_j) & \approx \frac{1}{\Gamma(1 - \gamma)}\int^{t_{j-1}}_0e^{-b(t_j - s)}\sum^{N_{exp}}_{k=1}
\omega_ke^{-s_k(t_j - s)}u'(s)ds\\
&\triangleq \frac{1}{\Gamma(1 - \gamma)}\sum^{N_{exp}}_{k=1}\int^{t_{j-1}}_0 \omega_ke^{-\tilde{s}_k(t_j -
s)}u'(s)ds~\left[\triangleq\frac{1}{\Gamma(1 - \gamma)}\sum^{N_{exp}}_{k=1}\omega_kU_{hist,k}(t_j)\right]\\
& = \frac{1}{\Gamma(1 - \gamma)}\sum^{N_{exp}}_{k = 1}\omega_k\left[e^{-\tilde{s}_k\tau}U_{hist,k}(
t_{j-1}) + \int^{t_{j-1}}_{t_{j-2}}e^{-\tilde{s}_k(t_j - s)}u'(s)ds\right],
\end{split}
\label{eqA.2}
\end{equation}
where $\tilde{s}_k = s_k + b$. To evaluate $U_{hist,k}(t_j)$ for $j = 1,2,\cdots,N_{exp}$, it
observes the following simple recurrence relation:
\begin{equation}
\begin{split}
U_{hist,k}(t_j) & = e^{-\tilde{s}_j\tau}U_{hist,k}(t_{j - 1}) + \int^{t_{j-1}}_{t_{j - 2}}
e^{-\tilde{s}_k(t_j - s)}u'(s)ds\\
& \approx e^{-\tilde{s}_k\tau}U_{hist,k}(t_{j-1}) + \frac{u(t_{j-1}) - u(t_{j-2})}{\tau}
\int^{t_{j-1}}_{t_{j-2}}e^{-\tilde{s}_k(t_j - s)}ds \\
& = e^{-\tilde{s}_k\tau}U_{hist,k}(t_{j-1}) + \frac{[u(t_{j-1}) - u(t_{j-2})](1 -
e^{-\tilde{s}_k\tau})}{\tau\tilde{s}_ke^{\tilde{s}_k\tau}}.
\end{split}
\label{eqA.3}
\end{equation}
Noting that $U_{hist,k}(t_1)\equiv 0$ when $n = 1$, we have
\begin{equation*}
{}^{FC}_{0}\mathbb{D}^{\gamma,\lambda(t)}_tu^{1} = \frac{u(t_1) - u(t_0)}{\tau\Gamma(2 - \gamma)}\Big(e^{-b\tau}\tau^{1 - \gamma} +
b\int^{\tau}_0e^{-b\theta}\theta^{1-\gamma}d\theta\Big),
\end{equation*}
where we define
\begin{equation}
{}^{FC}_{0}\mathbb{D}^{\gamma,\lambda(t)}_tu^{j} = \frac{u(t_j) - u(t_{j - 1})}{\tau\Gamma(2 - \gamma)}\Big(e^{-b\tau}\tau^{1 - \gamma} +
b\int^{\tau}_0e^{-b\theta}\theta^{1-\gamma}d\theta\Big) + \frac{1}{\Gamma(1 - \gamma)}
\sum^{N_{exp}}_{k=1}\omega_kU_{hist,k}(t_j)
\vspace{-4mm}
\end{equation}
as the approximate discrete operator for evaluating ${}^{C}_0D^{\gamma,\lambda(t)}_tu(t)\mid_{
t = t_j}$ quickly and $U_{hist,k}(t_j)$ can be computed via Eq. (\ref{eqA.3}). At each time step, we only need $\mathcal{O}(1)$ work to compute $U_{hist,k}(t_j)$ since $U_{hist,k}(t_{j-1})$ is known at that point. Thus, the total work
is reduced from $\mathcal{O}(M^2)$ to $\mathcal{O}(MN_{exp})$, and the total memory requirement is reduced
from $O(M)$ to $O(N_{exp})$\footnote{In our experiments, it always finds that $N_{exp} < 80$.}.

Similar to \cite{Jiang17}, replacing the \textit{L}1-type approximation (cf. Lemma \ref{lem2.1}) for the
generalized Caputo fractional derivative by our fast evaluation scheme ${}^{FC}_{0}\mathbb{D}^{\gamma,
\lambda(t)}_t$, we obtain a novel implicit difference scheme of the following form
\begin{equation}
\begin{cases}
{}^{FC}_{0}\mathbb{D}^{\gamma,\lambda(t)}_tu^{j+1}_i = \xi^{j+1}_i(\delta^{\alpha}_hu^{j+1}_i) + f^{j+1}_i,
& i = 1,2,\cdots,N-1,\quad j = 0,1,\cdots, M-1,\\
u^{0}_i = \phi(x_i),& i = 0,1,\cdots,N,\\
u^{j+1}_0 = \varphi(t_{j+1}),\quad u^{j+1}_N = \psi(t_{j+1}),& j = 0,1,\cdots,M-1,
\end{cases}
\label{eqA.5}
\end{equation}
which nearly reaches the approximation order of $\mathcal{O}(\tau^{2-\gamma} + h^2)$; see numerical results
in the next context. At each time step $t_{j+1}$, evaluating the right hand side (i.e., the known solutions at the
previous time levels) and inverting the linear system have $\mathcal{O}(NN_{exp})$ and $\mathcal{O}(I_{avg}
N\log N)$ computational complexity, respectively, which leads to an overall computational complexity of $
\mathcal{O}(MN(N_{exp} + I_{avg}\log N))$, where $I_{avg}(\ll N)$ is the average number of iterations required
for solving the resulting linear system at each time step. By contrast, if we use the Gaussian elimination method to solve the resulting linear systems of Eq. (\ref{eq1.9a}), the overall computational complexity of the implicit difference scheme (\ref{eq2.3}) is about $\mathcal{O}(MN^3 + M^2N)$ operations. In addition, it is meaningful to note that the above fast difference scheme has an overwhelming advantage when the number of temporal discretization steps (i.e., $M$) is relatively large.
\vspace{1mm}

\noindent\textbf{Example A.1} In this example, we test the fast difference scheme
(\ref{eqA.5}) and the direct difference scheme (\ref{eq2.3}) for solving the same
model problem in \textbf{Example 1} except different diffusion coefficient
$\xi(x,t) = 10(1/2 + x^2 + \sin t)$. Let the tolerance error $\epsilon
= 10^{-9}$ for fast difference scheme (\ref{eqA.5}) and Tables \ref{tabA.1}--\ref{tabA.2}
are reported to evaluate the accuracy and efficiency of the proposed algorithms.

\begin{sidewaystable}[!htpb]
\caption{$L_2$-norm and maximum norm errors versus grid size reduction when $\tau = 2^{-11}$ and $p = 0.7$ in Example A.1.}
\centering
\begin{tabular}{crcccccccccc}
\hline &\multicolumn{5}{c}{Direct scheme (\ref{eq2.3})} &\multicolumn{5}{c}{Fast scheme (\ref{eqA.5})}\\
[-2pt]\cmidrule(r{0.7em}r{0.7em}){3-7}\cmidrule(r{0.7em}r{0.6em}){8-12}\\[-11pt]
$(\gamma,\alpha,b)$ &$N$ &Error$_{\infty}$ &Rate$_{\infty}$ &
Error$_2$ & Rate$_2$ &CPU(s) &Error$_{\infty}$ &Rate$_{\infty}$ &Error$_2$ & Rate$_2$   &CPU(s)  \\
\hline
(0.2,1.1,1.0) &10 &7.3589e-2 &	--	 &7.0444e-2	&  --   &3.578 &7.3581e-2 &	--	  &7.0438e-2 & --	 &0.335 \\
              &20 &1.7410e-2 &2.0796 &1.7101e-2	&2.0424	&3.913 &1.7404e-2 &2.0799 &1.7095e-2 &2.0428 &0.408 \\
              &40 &4.1567e-3 &2.0664 &4.1035e-3	&2.0592	&4.159 &4.1515e-3 &2.0677 &4.0983e-3 &2.0605 &0.469 \\
              &80 &1.1354e-3 &1.8722 &9.8777e-4	&2.0546	&5.910 &1.1443e-3 &1.8592 &9.8257e-4 &2.0604 &1.157 \\
\hline
(0.5,1.5,1.0) &10 &4.8279e-2 &	--   &4.6207e-2	&  --   &3.586 &4.8274e-2 & --    &4.6202e-2 & --    &0.343 \\
              &20 &1.1381e-2 &2.0848 &1.0791e-2	&2.0983	&3.967 &1.1377e-2 &2.0851 &1.0787e-2 &2.0987 &0.417 \\
              &40 &2.7033e-3 &2.0738 &2.5503e-3	&2.0811	&4.278 &2.6990e-3 &2.0756 &2.5463e-3 &2.0828 &0.482 \\
              &80 &6.7243e-4 &2.0073 &6.0900e-4	&2.0662	&5.935 &6.7315e-4 &2.0034 &6.0516e-4 &2.0730 &1.161 \\
\hline
(0.9,1.9,1.0) &10 &4.6595e-2 &	--	 &4.6972e-2	&  --   &3.584 &4.6593e-2 & --    &4.6969e-2 & --    &0.328 \\
              &20 &1.1365e-2 &2.0356 &1.1402e-2	&2.0425	&3.966 &1.1363e-2 &2.0358 &1.1400e-2 &2.0427 &0.396 \\
              &40 &2.7816e-3 &2.0306 &2.7753e-3	&2.0386	&4.197 &2.7797e-3 &2.0313 &2.7734e-3 &2.0393 &0.473 \\
              &80 &6.8243e-4 &2.0272 &6.7697e-4	&2.0355	&5.857 &6.8051e-4 &2.0302 &6.7506e-4 &2.0386 &1.149 \\
\hline
(0.2,1.1,2.0) &10 &6.9685e-2 &  --	 &6.6419e-2	&  --   &3.285 &6.9679e-2 & --    &6.6415e-2 & --    &0.324 \\
              &20 &1.6438e-2 &2.0838 &1.6132e-2	&2.0417	&3.652 &1.6434e-2 &2.0840 &1.6128e-2 &2.0419 &0.373 \\
              &40 &3.9206e-3 &2.0679 &3.8721e-3	&2.0587	&3.983 &3.9167e-3 &2.0690 &3.8682e-3 &2.0598 &0.461 \\
              &80 &1.0747e-3 &1.8671 &9.3218e-4	&2.0544	&5.597 &1.0747e-3 &1.8657 &9.2820e-4 &2.0592 &1.112 \\
\hline
(0.5,1.5,2.0) &10 &4.5268e-2 &  --	 &4.3322e-2	&  --   &3.265 &4.5265e-2 &  --	  &4.3318e-2 &  --	 &0.321 \\
              &20 &1.0665e-2 &2.0856 &1.0114e-2	&2.0987	&3.512 &1.0662e-2 &2.0859 &1.0110e-2 &2.0992 &0.369 \\
              &40 &2.5332e-3 &2.0739 &2.3895e-3	&2.0816	&3.898 &2.5300e-3 &2.0753 &2.3865e-3 &2.0828 &0.458 \\
              &80 &6.2593e-4 &2.0169 &5.7014e-4	&2.0673	&5.697 &6.2648e-4 &2.0138 &5.6752e-4 &2.0722 &1.116 \\
\hline
(0.9,1.9,2.0) &10 &4.3556e-2 &  --   &4.3908e-2	&  --	&3.298 &4.3554e-2 &  --	  &4.3906e-2 &  --	 &0.319 \\
              &20 &1.0623e-2 &2.0357 &1.0657e-2	&2.0427	&3.557 &1.0621e-2 &2.0359 &1.0656e-2 &2.0428 &0.372 \\
              &40 &2.5992e-3 &2.0311 &2.5935e-3	&2.0388	&3.935 &2.5978e-3 &2.0316 &2.5921e-3 &2.0395 &0.463 \\
              &80 &6.3724e-4 &2.0282 &6.3217e-4	&2.0365	&5.713 &6.3580e-4 &2.0306 &6.3074e-4 &2.0390 &1.126 \\
\bottomrule
\end{tabular}
\label{tabA.1}
\end{sidewaystable}
\begin{sidewaystable}[!htpb]
\caption{$L_2$-norm, maximum norm errors and CPU time (in seconds) versus grid size reduction
when $N = \lceil2M^{(2-\gamma)/2}\rceil$ and $p = 0.7$ in Example A.1.}
\centering
\begin{tabular}{crcccccccccc}
\hline &\multicolumn{5}{c}{Direct scheme (\ref{eq2.3})} &\multicolumn{5}{c}{Fast scheme (\ref{eqA.5})}\\
[-2pt]\cmidrule(r{0.7em}r{0.7em}){3-7}\cmidrule(r{0.7em}r{0.6em}){8-12}\\[-11pt]
$(\gamma,\alpha,b)$ &$M$ &Error$_{\infty}$ &Rate$_{\infty}$ &
Error$_2$ & Rate$_2$ &CPU(s) &Error$_{\infty}$ &Rate$_{\infty}$ &Error$_2$ & Rate$_2$   &CPU(s)  \\
\hline
(0.2,1.1,1.0) &$2^5$ &3.2759e-3	&  --   &3.2241e-3 &  --   &0.008 &3.7666e-3 &  --	 &2.9594e-3	&  --	&0.007 \\
              &$2^6$ &1.0351e-3	&1.6621	&8.9483e-4 &1.8492 &0.029 &1.2961e-3 &1.5391 &7.6594e-4	&1.9500	&0.031 \\
              &$2^7$ &3.1910e-4	&1.6977	&2.4837e-4 &1.8491 &0.264 &4.4060e-4 &1.5566 &2.0299e-4	&1.9158	&0.254 \\
              &$2^8$ &9.5866e-5	&1.7349	&6.8823e-5 &1.8515 &2.873 &1.5160e-4 &1.5392 &5.8937e-5	&1.7842	&2.859 \\
\hline
(0.5,1.5,1.0) &$2^6$ &2.1236e-3	&  --	&2.0025e-3 &  --   &0.012 &2.0134e-3 &  --	 &1.9017e-3	&  --	&0.013 \\
              &$2^7$ &7.4051e-4	&1.5199	&6.7840e-4 &1.5616 &0.052 &7.5084e-4 &1.4231 &6.2573e-4	&1.6037	&0.053 \\
              &$2^8$ &2.7349e-4	&1.4370	&2.3228e-4 &1.5463 &0.356 &2.7741e-4 &1.4365 &2.0565e-4	&1.6053	&0.334 \\
              &$2^9$ &9.9638e-5	&1.4567	&8.0300e-5 &1.5324 &2.703 &1.0111e-4 &1.4561 &6.7254e-5	&1.6125	&2.442 \\
\hline
(0.9,1.9,1.0) &$2^7$    &5.7500e-3	&  --	&5.7539e-3 &  --   &0.019 &5.7304e-3 &  --   &5.7345e-3	&  --	&0.021 \\
              &$2^8$    &2.5255e-3	&1.1870	&2.5185e-3 &1.1920 &0.078 &2.5143e-3 &1.1885 &2.5073e-3	&1.1935	&0.054 \\
              &$2^9$    &1.1846e-3	&1.0922	&1.1776e-3 &1.0967 &0.303 &1.1782e-3 &1.0936 &1.1713e-3	&1.0980	&0.161 \\
              &$2^{10}$ &5.3828e-4	&1.1380	&5.3354e-4 &1.1422 &1.776 &5.3475e-4 &1.1396 &5.3004e-4	&1.1439	&0.705 \\
\hline
(0.2,1.1,2.0) &$2^5$ &3.1058e-3	&  --	&3.0395e-3 &  --   &0.007 &3.3736e-3 &  -- 	 &2.8288e-3 &  --	&0.008  \\
              &$2^6$ &9.7378e-4	&1.6733	&8.4355e-4 &1.8493 &0.028 &1.1535e-3 &1.5483 &7.4229e-4	&1.9301	&0.032  \\
              &$2^7$ &2.9435e-4	&1.7261	&2.3410e-4 &1.8493 &0.261 &3.8784e-4 &1.5725 &2.0336e-4	&1.8679	&0.255  \\
              &$2^8$ &8.8505e-5	&1.7337	&6.4856e-5 &1.8518 &2.901 &1.3145e-4 &1.5609 &5.5811e-5	&1.8654	&2.861  \\
\hline
(0.5,1.5,2.0) &$2^6$ &1.9855e-3	&  --   &1.8722e-3 &  --   &0.011 &1.9019e-3 &  --   &1.7954e-3	&  --	&0.012 \\
              &$2^7$ &6.8955e-4	&1.5258	&6.3406e-4 &1.5620 &0.053 &6.9740e-4 &1.4474 &5.9410e-4	&1.5955	&0.054 \\
              &$2^8$ &2.5469e-4	&1.4369	&2.1704e-4 &1.5467 &0.360 &2.5766e-4 &1.4365 &1.9681e-4	&1.5939	&0.255 \\
              &$2^9$ &9.2793e-5	&1.4567	&7.5014e-5 &1.5327 &2.699 &9.3904e-5 &1.4562 &6.6809e-5	&1.5587	&2.434 \\
\hline
(0.9,1.9,2.0) &$2^7$    &5.3618e-3 &  --   &5.3656e-3 &  --	  &0.020 &5.3470e-3	&  --   &5.3510e-3 &  --   &0.021 \\
              &$2^8$    &2.3546e-3 &1.1872 &2.3482e-3 &1.1922 &0.077 &2.3462e-3	&1.1884	&2.3398e-3 &1.1934 &0.055 \\
              &$2^9$    &1.1044e-3 &1.0922 &1.0980e-3 &1.0967 &0.302 &1.0996e-3	&1.0933	&1.0933e-3 &1.0977 &0.255 \\
              &$2^{10}$ &5.0182e-4 &1.1380 &4.9742e-4 &1.1423 &1.769 &4.9917e-4	&1.1394	&4.9480e-4 &1.1438 &0.694 \\
\bottomrule
\end{tabular}
\label{tabA.2}
\end{sidewaystable}
Tables \ref{tabA.1}--\ref{tabA.2} illustrate the temporal/spatial errors, convergence orders and CPU time of the
methods. It can be seen from Table \ref{tabA.1} that when $\tau = 2^{-11}$, both ``${\rm Error}_{\infty}$"
and ``${\rm Error}_2$" of two implicit difference schemes for the variable coefficient GTSFDEs with different $(\gamma,
\alpha,b)$'s decreases steadily for smaller $h$, and the order of accuracy in space is about two. Fixing $N = \lceil
2M^{(2-\gamma)/2}\rceil$, Table \ref{tabA.2} lists the maximum-norm and $L_2$-norm errors and illustrates that the
order of temporal accuracy is of $(2 - \gamma)$. Therefore, Tables \ref{tabA.1}--\ref{tabA.2} confirm that
the rate of the truncation errors of numerical schemes (\ref{eq2.3}) and (\ref{eqA.5}) is $\mathcal{O}(\tau^{2-\gamma}
+ h^2)$. However, it seems that the temporal errors based on the maximum norm of fast scheme (\ref{eqA.5}) change slightly irregularly compared
to those of the direct scheme (\ref{eq2.3}), especially for the case of $(0.2,1.1,2.0)$.
Moreover, the fast scheme
(\ref{eqA.5}) requires less CPU time than the direct scheme (\ref{eq2.3}) for the variable-coefficient GTSFDEs with
different $(\gamma,\alpha,b)$'s. The time reduction between the direct scheme (\ref{eq2.3})
and the fast scheme (\ref{eqA.5}) shown in Table \ref{tabA.2} is not distinct, because the
number of temporal discretization steps is less than the size of spatially discretized linear systems which are preponderantly time-consuming. In conclusion, although the derived
fast scheme (\ref{eqA.5})\footnote{In fact, the above fast scheme can easily utilize the non-uniform temporal steps \cite{Stynes17,Liao18}, which can enhance its (numerical)
temporal convergence order for solving the variable-coefficient GTSFDEs (even with the weak singularity at initial time). MATLAB codes of all the numerical tests is available from the authors' emails.} needs less CPU time and memory cost than the direct scheme (\ref{eq2.3}). Further analysis is still required to assess its stability and convergence properties.
\section*{Acknowledgments}
{\em The authors are grateful to Prof. Jiwei Zhang and Dr. Hong-Lin Liao for their
constructive discussions and insightful comments. This research is supported by NSFC
(11801463 and 61772003), the Applied Basic Research Project of Sichuan Province (20YYJC3482),
the Fundamental Research Funds for the Central Universities (JBK1902028)
and the Ministry of Education of Humanities and Social Science Layout Project (19JYA790094).
Meanwhile, the first author would like to thank Prof. Hai-Wei Sun and Dr. Siu-Long
Lei for their helpful discussions during his visiting to the University of Macau.
}
\medskip
\medskip

\noindent\textbf{References}

\end{document}